\documentclass[12pt,reqno]{amsart}
\usepackage[a4paper,top=2.8cm,bottom=3cm,left=2.7cm,right=2.7cm,marginparwidth=55pt]{geometry}

\usepackage{latexsym}
\usepackage{amsthm, amsmath, amssymb,latexsym,amsfonts,verbatim}
\usepackage{comment}

\usepackage{tikz-cd}
\usetikzlibrary{cd,shapes,arrows,positioning,calc}

\usepackage{paralist}
\usepackage[mathscr]{eucal}

\usepackage[colorinlistoftodos]{todonotes}
\usepackage[all]{xy}

\newcommand{\sym}{\mbox{Sym}}
\renewcommand{\H}{\mbox{H}}
\newcommand{\surjarrow}{\rightarrow\!\!\!\!\rightarrow}

\usepackage[colorlinks,bookmarks,pdftex]{hyperref} %
\usepackage[capitalise]{cleveref}
\xyoption{all} \numberwithin{equation}{section}
\newtheorem{theorem}[equation]{Theorem}
\newtheorem{corollary}[equation]{Corollary}
\newtheorem{proposition}[equation]{Proposition}
\newtheorem{lemma}[equation]{Lemma}

\theoremstyle{definition}

\newtheorem{definition}[equation]{Definition}

\newtheorem{remark}[equation]{Remark}

\newtheorem{setting}[equation]{Setting} 

\newtheorem{problem}[equation]{Problem}

\newif\ifpdf\ifx\pdfoutput\undefined\pdffalse\else\pdfoutput=1\pdftrue\fi

\newcommand{\nc}{\newcommand} 
\nc{\cH}{{\mathcal H}}
\nc{\cA}{{\mathcal A}}
\nc{\cG}{{\mathcal G}}
\nc{\cC}{{\mathcal C}}
\nc{\cO}{{\mathcal O}}
\nc{\cI}{{\mathcal I}}
\nc{\cR}{{\mathcal R}}
\nc{\cB}{{\mathcal B}}
\nc{\cY}{{\mathcal Y}}
\nc{\cK}{{\mathcal K}} 
\nc{\cX}{{\mathcal X}}
\nc{\cS}{{\mathcal S}}
\nc{\cE}{{\mathcal E}}
\nc{\cF}{{\mathcal F}}
\nc{\cZ}{{\mathcal Z}}
\nc{\cQ}{{\mathcal Q}}
\nc{\cN}{{\mathcal N}}
\nc{\cP}{{\mathcal P}}
\nc{\cL}{{\mathcal L}}
\nc{\cM}{{\mathcal M}}
\nc{\cT}{{\mathcal T}}
\nc{\cW}{{\mathcal W}}
\nc{\cU}{{\mathcal U}}
\nc{\cJ}{{\mathcal J}}
\nc{\cV}{{\mathcal V}}
\nc{\bH}{{\mathbb H}}
\nc{\bA}{{\mathbb A}}
\nc{\bG}{{\mathbb G}}
\nc{\bC}{{\mathbb C}}
\nc{\bO}{{\mathbb O}}
\nc{\bI}{{\mathbb I}}
\nc{\bB}{{\mathbb B}}
\nc{\bY}{{\mathbb Y}}
\nc{\bK}{{\mathbb K}} 
\nc{\bX}{{\mathbb X}}
\nc{\bS}{{\mathbb S}}
\nc{\bE}{{\mathbb E}}
\nc{\bF}{{\mathbb F}}
\nc{\bZ}{{\mathbb Z}}
\nc{\bQ}{{\mathbb Q}}
\nc{\bN}{{\mathbb N}}
\nc{\bP}{{\mathbb P}}
\nc{\bL}{{\mathbb L}}
\nc{\bM}{{\mathbb M}}
\nc{\bT}{{\mathbb T}}
\nc{\bW}{{\mathbb W}}
\nc{\bU}{{\mathbb U}}
\nc{\bD}{{\mathbb D}}
\nc{\bJ}{{\mathbb J}}
\nc{\bV}{{\mathbb V}}
\nc{\bbZ}{{\mathbb Z}}
\nc{\bR}{{\mathbb R}}
\nc{\fr}{{\rightarrow}}
\nc{\co}{{\nabla}}
\renewcommand{\H}{H}
\newcommand{\la}{\longrightarrow}
\nc{\cu}{{\barline{\nabla}}}
\newcommand{\pr}{\mathbb P}
\newcommand{\pic}{\mbox{\upshape{Pic}}}
\nc{\divi}{\mbox{\upshape{div}}}

\newcommand{\derx}[1]{\frac{\partial #1}{\partial x}}

\begin{document}



\title {Conic linear series and pencils of plane quartics} 
\author{R.~Moschetti, G.P.~Pirola, L.~Stoppino}

\begin{abstract}
    We study linear systems cut out by cones on a smooth complex non-degenerate curve $C$ in $\bP^{3}$. 
    We develop a systematic study of families of such systems, considering their limits, their infinitesimal behaviour, and some associated geometric structures. As an application, by using the machinery developed, we prove the existence of a non-isotrivial pencil of quartics with only one base point, all whose members are irreducible and whose general member is smooth.
\end{abstract}
\maketitle
{\footnotesize 2010 Mathematics Subject Classification 14H50; 14C21; 14H10; 14C20.}

{\hypersetup{linkcolor=black}\tableofcontents}

\section{Introduction}\label{intro}
Let $C\subset\bP^{3}$ be a smooth non-degenerate  complex projective curve of degree $d$. 
We call a \emph{conic linear system} on $C$ a linear system which is cut out on $C$ by cones of fixed degree whose vertex is disjoint from $C$ (see \Cref{limnormali}). 
In this paper we carry out a detailed study of these objects and of their limits as the vertex specializes to a point of $C$.

Our original motivation comes from the following problem: to construct, if possible, non-isotrivial pencils of plane curves with a single base point, whose members are all irreducible and whose general member is smooth. 
The case of cubics, which we studied extensively in the recent paper \cite{MPS},  appeared many  times in the literature from different perspectives: starting with Cayley \cite{Cayley} and ending - to our knowledge - with Koll\'ar \cite[Example 46]{Kollar}. 
The construction in \cite{MPS} is related to the torsion points of order $9$ on a plane cubic (see also \cite{Gat}). 

We wondered if the method of  \cite{MPS} could also be modified to treat the case of higher degree plane curves. 
In order to work with torsion points as before, it seemed unavoidable to use elliptic curves.  
We thus tried to work in $\pr^{n-1}$ with elliptic curves of degree $n$ and hypersurfaces, coming back to $\pr^2$ by projections. 
This naturally led us to study cones, and divisors cut out on the curve by them. 
We soon realized that our wishful thinking was difficult to implement in full generality. 
In this paper we set the theoretical framework to work with curves of any degree in $\pr^3$, obtaining some general results on cones and ``conic'' divisors for any smooth curve in $\pr^3$. 
As an application we solve our original problem for $d=4$.

\smallskip

We now describe in detail the content of the paper. 
As a first result, we describe linear limits of cones over a curve $C$.  
This result will be used extensively later in the paper when working with divisors.
\begin{proposition}\label{prop1} [\Cref{lemma:limitconesexplicit}]
Let $p\in C \subset \pr^3$ and let $t_p$ be the projective tangent line to $C$ at $p$. 
Consider a line $\ell\neq t_p$ through $p$. Then the limit of the degree-$d$ cone over $C$, as its vertex specializes to $p$ along $\ell$, is the union of the cone over $C$ with vertex $p$ and the plane spanned by $\ell$ and $t_p$.
\end{proposition} 
This result is closely related to \cite[Section 3]{CMR}, see \Cref{rem: ciliberto}.
We give two proofs of \Cref{prop1}: one geometric in \Cref{lemma:limitconesexplicit} and one analytical in \Cref{subs:coordinates}.
While the second approach is less elegant, it introduces techniques that are essential in the rest of the paper.

\smallskip

Let  $\nu\colon\widetilde \pr^3\to \pr^3$, the blow-up of $\pr^3$ along $C$.  
One of the outcomes of \Cref{lemma:limitconesexplicit} is a natural injective morphism from $\widetilde \pr^3$ to a projective space.
Let us explain briefly
this construction, described in \Cref{sec: proj model}. 
Denote by $U\subseteq \pr^3\smallsetminus C$ the locus of points from which the linear projection map restricted to $C$ is birational. 
Let $I$ be the ideal of $C$ in $\pr^3$ and let $c\colon U \to \pr(I(d))$ be the map assigning to each $p \in U$ the cone over $C$ with vertex at $p$. 
Let $E=\pr(N_{C|\bP^3})$ be the exceptional divisor of $\nu$. Any $\xi\in E$ corresponds to a point $p\in C$ and to a plane $H$ containing $t_p$. Proposition \ref{prop1} gives the extension of $c$ along $E$, while the limits at the finite set of non-birational projection points outside $C$ are unique.
 
In \Cref{extension c} we prove that $c$ extends to an injective morphism $\widetilde c\colon\widetilde \pr^3\to\pr(I(d))$, and we compute the associated divisor class in $\pic(\widetilde \pr^3)$.

\smallskip

We then turn our attention to the structure of {\em conic linear systems}:
denote by $R_k(q)$ the linear system cut on $C$ by the cones of degree $k$ in $\pr^3$ with vertex $q\in U$, see \Cref{limnormali}.

For every integer $1\leq k\leq d$, set  $n_k=\dim R_k(q)$ and define the {\em cone map}
\[
\rho_k\colon U\longrightarrow G(n_k,H^0(C,\cO_C(k)))
\]
where $\rho_k(q):=R_k(q)$. 

We are interested in defining a limit system $R_k^\xi(p)$ for $\xi=(p,H)\in E$, or equivalently in extending $\rho_k$ to $\widetilde \pr^3$. 
Let $C'$ be the locus of points of $p\in C$ such that  the projection from $p$ is birational onto its image. 
For $k\leq d-2$ and $p\in C'$, the linear limit along any line $\ell\neq t_p$ is simply $R_k(p)$.

On the other hand, for $k\geq d-1$ there is not an easy construction as the one made for the map $c$. In particular, we can define linear limits of conic series  $R^\ell_k(p)$ along a line $\ell\neq t_p$ passing through $p$, but it is not clear that these limits only depend on the 
point $\xi=(p,H)\in E$, where $H=\langle \ell, t_p\rangle$. 
Indeed, we can not prove that this definition is everywhere independent on the choice of $\ell$, but we can give some partial results for the cases $k=d-1$ and $k=d$.

For $k=d-1,d$ let $\widetilde \rho_k$ be the maximal  extension of $\rho_k$, and let  $\cV_k$ the maximal open subset of $\widetilde{\pr^3}$ where it is defined:
\begin{equation}\label{rho-intro} \widetilde{\rho_k}\colon \cV_k\longrightarrow    G(n_k,H^0(C,\cO_C(k))).\end{equation}
Consider the finite set $\Gamma:=(C\smallsetminus C')\cup 
\{p\in C' \mid t_p\cap C\neq \{p\}\}\subset C.$ Call $\mu\colon E\longrightarrow C$ the natural projection.
In \Cref{sec: limits} we prove the following. 
\begin{theorem}[\Cref{teo: VsuE}]
We have that 
$\cV_k\cap E$ is contained in the pullback $\mu^*(\Gamma)$.
\end{theorem}
Given $\xi \in  \cV_k\cap E$, we prove  that $R_k(\xi)=R^\ell_k(p)$, where $\xi=(p,H)$ and $H=\langle \ell, t_p\rangle$.

Moreover, we show in \Cref{cor:evaluation-d-1} and  \Cref{limd} that, every element of $\pr(R^\xi_k(p))$ has evaluation at least $d-2$ at $p$. 
In our arguments the analytic approach started in the first section is crucial. 

\smallskip

In \Cref{sec: conemap} we study the  differential of $\rho_d$, proving in particular:
\begin{theorem}[\Cref{limb}]
The differential of $\rho_d$ is injective over $U$.
 \end{theorem}

Eventually, in \Cref{sec: conic divisors} we encode the geometrical structure of the conic linear series in a morphism
$$\Phi: \bP_U(\cE)\la \bP(\H^0(C,\cO_C(d))),$$
where $\cE$ is a vector bundle whose fibres are linear systems of degree $d$ cones with vertex  $p\in U$, modulo those cones that contain the curve $C$.
So, the fibres are naturally isomorphic to $ R_d(p)$.
The morphism $\Phi$ assigns to each class of cones the conic divisor cut out on $C$.
By using a combination of classical methods (e.g. Castelnuovo's bound in \Cref{d(imp)}) and some subtle computations we study the differential of $\Phi$:
\begin{theorem}[\Cref{thm:diffmaxrank}]
    The differential of $\Phi$ is generically of maximal rank. For a general $y \in \bP_U(\cE)$ the following holds:
    \begin{enumerate}
\item if $d\geq 5$ or $d=4$ and $g=0,$ then $d\Phi_y$ is generically injective;
\item  if $d=4$ and $g=1$, then $\dim \ker d\Phi_y =1;$
\item  if $d=3$ and $g=0$, then $\dim \ker d\Phi_y=2$.
\end{enumerate}
\end{theorem}

This yields the following result:
\begin{corollary}[\Cref{cor: dominant}]
The following facts are equivalent:
\begin{enumerate}
    \item[(i)] the map $\Phi$ is dominant,
    \item[(ii)] the degree $d\leq 4$.
\end{enumerate}
\end{corollary}

Using the map $\widetilde{\rho_d}$ of \Cref{rho-intro}, we obtain a partial compactification of the construction and define a partial extension $\Phi'$ of $\Phi$, whose domain parametrizes linear limits of systems of cones. 
We then extend further $\Phi$ to a morphism $\Phi''$ from  $\pr_{\widehat{\pr^3}}$, where $\widehat{\pr^3}$ is a suitable blow up of $\widetilde {\pr^3}$ (\Cref{def: phi''}).

\smallskip 

In \Cref{sec:application} we turn our attention to the construction of pencils as explained in the beginning of this introduction. 
We consider the case where  $C\subset \pr^3$ is an  elliptic smooth curve of degree $4$.
In this case we can prove that the map $\Phi''$ is indeed surjective (\Cref{do}). 
By exploiting this  result together with the geometric structure of these maps, we show the following:
\begin{theorem}\label{teo1}(\Cref{teo: contro})
Let $C\subset \pr^3$ be a smooth elliptic curve of degree $4$ and let $D$ be any effective divisor $D\in|\mathcal O_C(4)|$.
There is a $1$-dimensional family of cones in $B\subset \bP_U(\cE)$ 
such that:
\begin{enumerate}
\item[(i)] the fibre projection of $B$ in $U$ is again $1$-dimensional;  
\item[(ii)] by calling $\widetilde{K}_t$ the cone corresponding to $t\in B$, we have that  $\Phi(\widetilde{K}_t)= D$, for any $t\in B$.
\end{enumerate}
\end{theorem}

Finally,  we obtain our geometric application:
\begin{theorem}[\Cref{thm: pencil}]\label{thm: 1}
There exists a pencil of quartics in $\bP^2$ such that:
 \begin{enumerate} 
  \item the base locus is set-theoretically reduced to one point;
 \item all the quartics of the pencil are irreducible;
\item the general element of the pencil is smooth;
 \item the pencil is non-isotrivial.
  \end{enumerate}
 \end{theorem}
The idea of the proof of this last result is the following. 
Let $C\subset \pr^3$ be a smooth elliptic quartic. 
Fix the origin of its group structure at a flex point $O$. Let $\sim$ denote the linear equivalence between divisors on $C$.
Choose a point $q\in C$ such that $16 q\sim 16O$ but $8q\not\sim 8O$ (in direct analogy with the construction used for plane cubics in \cite{MPS}). 
By \Cref{teo1} there exists a $1$-dimensional family of cones in the inverse image of $16q$ under the map $\Phi$. 
Fixing one such cone and projecting from its vertex we obtain a plane quartic which intersects the projection of $C$ precisely on the image of $q$. 
These two curves generate the desired pencil. 
The irreducibility of all its members comes from the assumption we made on $q$. 
The smoothness of the general member relies on a delicate argument that makes heavy use of all the results of the previous sections. 
In particular, we will need the aforementioned result stating that every element of $\pr(R^\xi_4(p))$ has evaluation at $p$ at least $2$ (\Cref{limd}), the existence of the $1$-parameter family of cones proved in \Cref{teo1}, and an explicit description of the inverse image ${\Phi''}^{-1}(16q)$ given in \Cref{lemma: conovero}.
The non-isotriviality is ensured by the existence of a reduced fibre with geometric genus $1$ (the image of $C$).

\smallskip

We believe that conic linear systems are interesting by themselves and will have further applications beyond the one given in this work. 
In the last section we list some open questions and describe some possible developments that our theory can have.

\subsection*{Acknowledgements}
We thank Letterio Gatto and Joan Carlos Naranjo for useful conversations on the subject of this paper. We are grateful to  Ciro Ciliberto for invaluable insight on this work, in  particular for suggesting some possible developments (\Cref{cir1}, \Cref{cir2}). 
We also wish to thank  Wei Chen for careful reading of the paper and for pointing out an inaccuracy in a first version, that is mended here in \Cref{sec: limits}.
All the authors are members of the GNSAGA - INdAM and UMI and are partially supported by the local research funding FAR UniPV.

\section{Preliminaries} \label{sec:preliminaries}
Let $C\subset \bP^3$ be a smooth irreducible non-degenerate curve of genus $g$ and degree $d$ defined over the field $\bC$ of complex numbers. Consider the short exact sequence 
$$0\longrightarrow I\longrightarrow\cO_{\bP^3}\longrightarrow\cO_C\longrightarrow 0,$$
where $I$ is the ideal sheaf of $C$. Let $\cO_{\bP^3}(1)$ denote the tautological line bundle and let $\cO_C(1)$ be its restriction to $C$. Let $V\subset \H^0(C,\cO_C(1))$ be the $4$-dimensional vector space that defines the embedding $C \subset \bP^3 = \bP(V^\ast)$. 
In particular, $|V|=\bP(V)\cong \bP^3$ is the linear system on $C$ of divisors given by intersecting $C$ with the planes of $\bP^3$. 
A point $p\in \bP(V^\ast)$ corresponds to a $3$-dimensional subspace $W(p)\subset V$, namely $W(p)$ consists of the sections of $V$ vanishing at $p$ (see also \Cref{defn:VmD}). 
Dually, we have a surjection $V^\ast\surjarrow W(p)^\ast$ and a projection map 
\[\Pi_p\colon \pr(V^*)\dasharrow \pr(W(p)^\ast),\]
defined on $ \pr(V^\ast)\smallsetminus \{p\}$.
Consider the restriction to $C$ of $\Pi_p$, which we denote by $\pi_p\colon C \dasharrow \bP(W(p)^\ast)$. 
When $p \notin C$, the map $\pi_p$ is already a morphism, while when $p \in C$, as $p$ is a smooth point of $C$, $\pi_p$ can be extended to a morphism. 
We will call  $C_p:=\pi_p(C)\subset \pr(W(p)^\ast)$ the image of this morphism.

\smallskip

It will be useful to consider a partition of $\bP^3$ into three subsets, according to the properties of the map $\pi_p$.
Let us call $S\subset \pr^3$  the set of points in which the morphism $\pi_p$ is not birational, let  $C'$ be  the set $C \smallsetminus S$, and $U$ be the complement $U=\bP^3\smallsetminus (C \cup  S)$, so that we have
\begin{equation} \label{dec}
\bP^3=S \sqcup C' \sqcup U.
\end{equation}

\begin{lemma} \label{lem:suc}
We have the following:
\begin{enumerate}
\item $p\in U \iff \pi_p$ is birational and $C_p$ has degree $d$;
\item $p\in C'\iff \pi_p$ is birational and $C_p$ has degree $d-1$;
\item $p\in S \iff \pi_p$ is not birational. 
\end{enumerate}
Moreover, $S$ is a finite set of points.
\end{lemma}
\begin{proof}
The first three claims follow from the properties of linear projections. The finiteness of $S\cap C$ comes from the trisecant lemma \cite[Chap.II, Sec.5]{GH}. 
For the finiteness of $S\smallsetminus C$ we can use the result of 
\cite{PS}.
An alternative argument can be done by means of the theory of focal loci that is presented for instance in \cite{CF}, see also \cite{CCM, CM}.
\end{proof}

\begin{remark}\label{rem: mah}
For $p\in U$, the plane curve $C_p$ is singular and reduced of degree $d$,
and when $p\in C'$, the curve $C_p$ is reduced of  degree $d-1$. 
\end{remark}

For any integer $k\geq 1$ let us consider the space $V_k :=\sym^k(V)=\H^0(\bP^3,\cO_{\bP^3}(k))$, which parametrizes hypersurfaces of degree $k$ in $\bP^3$. 
For a linear subspace $W \subset V$, we have $W_k:=\sym^k(W) \subset V_k$,  parametrizing hypersurfaces of degree $k$ in $\bP(W^\ast)$, or, equivalently, {\em cones} of degree $k$ in $\bP(V^\ast)$ with vertex  $\bP Ann(W)\subset \bP(V^\ast)$. 

\begin{lemma}\label{lem: dimensioni}
We have the following characterization of the decomposition \eqref{dec}:
\begin{enumerate}
\item $p\in U \iff \dim W(p)_d\cap I(d)=1$;
\item $p\in C'\iff{\dim W(p)_d}\cap I(d)=3$;
\item $p\in S \iff \dim W(p)_d\cap I(d) \geq 6.$
\end{enumerate}
\end{lemma}
\begin{proof}
Let $p\in \bP^3$. Let $C_p\subset \pr(W(p)^\ast)\cong \bP^2$ be the image curve and $f_p=0$  an associated equation. 
As \eqref{dec} is a partition, it is enough to prove all the implications from left to right in the three points of the statement. 

$(1)$ If $p\in U$, then $f_p$ has degree $d$, and it represents the unique cone with vertex $p$ over $C$. So $W(p)_d\cap I(d)=\langle f_p\rangle$.

$(2)$ If $p\in C'$, then the cones with vertex $p$ over $C$ all have equation $f_ph=0$ where $h\in W(p)$. Hence its dimension is $3$.

$(3)$ If $p\in S$, then $C_p$ has degree $k<d-1$. So the space $W(p)_d\cap I(d)$ consists in the polynomials of the form $f_pg$, where $g\in W(p)_{d-k}$. 
Now we have, as wanted,
$$\dim W(p)_{d-k}=\binom{d-k+2}{2}\geq \binom{4}{2}=6.$$ 
\end{proof}

\begin{definition}\label{def: coni}
For any $p\in \pr^3$, we denote by  $\widetilde C_p\subset \bP^3$ the cone over $C$ with vertex $p$. Let  $0 \neq f_p$ be a polynomial describing the hypersurface $\widetilde C_p$. 
 \end{definition}
 The polynomial $f_p$ is defined up to a scalar and has degree $d$ if $p\in \pr^3\smallsetminus C$, degree $d-1$ for $p\in C$. When $p\in S$, the cone $\widetilde C_p$ is non-reduced; see \Cref{rem: mah}.
 
 \begin{remark} \label{rmk:fixedcone}
 Recalling that $\Pi_p\colon \pr(V^*)\dasharrow \pr(W(p)^\ast)$ is the projection map, 
we have that set-theoretically \[
\widetilde C_p=\overline{{\Pi_p}^{-1}(C_p).}
\]
By choosing a plane $\Sigma\subset \pr^3$ not containing $p$ in $\pr^3$ and identifying it with $ \pr(W(p)^\ast)$, one can also think as $\widetilde C_p$ as the cone in $\pr^3$ with base the plane curve $C_p\subset \Sigma$ and vertex $p$. 
Moreover, $f_p=0$ can be thought as both: (1) the equation of the cone with vertex $p$ over $C$ in $\bP^3$; (2) the equation of the projection curve $C_p \subset \bP(W(p)^\ast)\cong\bP^2$.
\end{remark}
\begin{definition}\label{def: c map}
If $p\in U$, the polynomial $f_p$ is a generator of  $W(p)_d \cap I(d)$, see \Cref{lem: dimensioni}.  
The correspondence $p\mapsto W(p)_d\cap I(d)$ thus defines a map
\begin{align} \label{c map}
\begin{split}
c:U &\longrightarrow \bP (I(d)) \\
 p & \,\mapsto \, \,[f_p].
 \end{split}
\end{align}
\end{definition}
In \Cref{extension c} we will extend $c$ to the blow-up of $\pr^3$ along $C$ and compute the associated linear system.

\section{Linear limits of cones} 
There are two families of dimension $3$ of cones of degree $d$ containing $C$: the first one consists of cones $\widetilde C_p$ with $p\in U$, the second one is made of cones of the form $\widetilde C_p\cup H$ with vertex in $p\in C$, and $H$ a plane belonging to $W(p)$. 
As a consequence the elements of the second family cannot  all be limits of elements of the first one. 
In this section we study the limits of the cones of the first family as the vertex specializes linearly  to a point of $C'$ and see that they become certain  elements of the second family. 

Fix a point $p\in C$. The limit cone for a point moving to $p$ depends on the direction. Let us consider the following situation: suppose we have a point $p_t\in \pr^3$ going to $p$ along a line $\ell \subset \bP^3$. Consider $\Delta$ a small disk in $\bC$ centred in $p$, and an affine coordinate $t\in \Delta$. The line $\ell$ near $p_t$ is parametrized by the affine coordinate $t$ as
\begin{equation} \label{eqn:retta}
\ell=\{ p+tv, t\in \bC\},
\end{equation}
where $v$ is the direction of $\ell$. Suppose that for any $t\in \Delta ^\ast$ we have $p_t\in U$. 
Then the polynomial $f_{t}:=f_{p_t}$ of \Cref{def: coni} belongs to $V_d$ for $t \neq 0$.

\begin{lemma} \label{lemma: limitcones}
Given $p\in C$ on a line $\ell \subset \bP^3$, there exists a unique (flat-)limit of $[f_{t}]$ in $\bP(V_d)$ for $t \to 0$, that is a degree $d$ cone with vertex $p$ and containing $C$.
\end{lemma}
\begin{proof}
First observe  that the limit of the family $[f_{t}]\in \bP(V_d)$ exists  by the valuative criterion of properness, and that the locus 
\[
\mathcal C:=\{[f]\in \bP(V_d)\mid Z(f) \mbox{ is a  cone} \}=\cup_{p\in \pr^3}W(p)_d
\]
is closed in $ \bP(V_d)$, so the limit also is a degree $d$ cone. 
Finally, it is clear that the vertex of this limit cone contains $p$, and that the limit cone contains $C$. 
\end{proof}

\begin{definition}\label{defn:limitcones}
We call the construction  above \emph{linear limit of $\widetilde C_{p_t}$}.
We denote the corresponding hypersurface as $\widetilde{C}^\ell_0$.
\end{definition}

As in \Cref{rmk:fixedcone}, we fix a general plane $\mathbb{P}^{2}\subset\mathbb{P}^3$ as the target space. We will use the same notation $\Pi_{p_t}\colon\bP^3\smallsetminus\{p_t\} \to \bP^{2}$.
The cone $\widetilde C_{p_t}$ is completely determined by the data $(p_t,\Pi_{p_t}(C))$; indeed, it is precisely the cone with vertex $p_t$ over $\Pi_{p_t}(C)$, and this remains true for $t=0$. 

\begin{proposition} \label{lemma:limitconesexplicit}
Let $p\in C$ and let $t_p$ be the projective tangent line to $C$ at $p$. 
Let $\ell \neq t_p$ be a line passing through $p$. Then
$$\widetilde{C}^\ell_0= \widetilde C_p\cup H,$$ 
where $H=\langle \ell,t_p \rangle$ is the unique plane in $\bP^3$ containing $t_p$ and $\ell$.
\end{proposition}
\begin{proof}
First, we prove the result in the special case when $\ell$ meets $C$ only in $P$.

Note first that the limit cone $\widetilde C_0^\ell$ contains the cone $\widetilde C_p$, which has degree $d-1$. As a consequence, we must have
$\widetilde C_0^\ell=\widetilde C_p\cup K$
for some plane $K$ through $p$.

Moreover, for every $t\neq 0$ the line $\ell$ is contained in $\widetilde C_{p_t}$, hence also in the limit $\widetilde C^\ell_0$. Therefore $K$ contains $\ell$. It remains to show that $K$ also contains $t_p$, i.e. that $K=\langle \ell,t_p\rangle$.

Let $\Lambda\simeq \mathbb P^2$ be a fixed plane disjoint from $p$, and consider the projections to $\Lambda$. Intersecting $\widetilde C^\ell_0$ with $\Lambda$, we get
$\widetilde C^\ell_0\cap \Lambda=\Pi_p(C)\cup K'$, where $K':=K\cap \Lambda$.
Here $\Pi_p(C)$ has degree $d-1$, while $K'$ is a line.

Set $q:=\ell\cap \Lambda$, which coincides with $\Pi_{p_t}(p)$ when $t \neq 0$.
For every $t\neq 0$, the projected curve $\Pi_{p_t}(C)$ is then tangent to the line $\Pi_{p_t}(t_p)$ at $q$. Notice that this line is independent of $t$: indeed, it is precisely the line $H\cap \Lambda$,
since it is the image of $t_p$ under projection from $p_t\in \ell$, and $H=\langle \ell,t_p\rangle$.

Now $q\notin \Pi_p(C)$. Indeed, the fibre of $\Pi_p$ over $q$ is exactly $\ell$, and by assumption $\ell\cap C=\{p\}$. The projection from $p$ is not defined at $p$, so the point $q$ does not belong to the image of $C$. More precisely, the closure of $\Pi_p(C)$ is obtained by adding the point $t_p\cap \Lambda$, which is different from $q$ because $\ell\neq t_p$.

Therefore, in the limit curve $\widetilde C^\ell_0\cap \Lambda=\Pi_p(C)\cup K'$,
the only component passing through $q$ is the line $K'$. Since the curves $\Pi_{p_t}(C)$ are tangent at $q$ to the fixed line $H\cap \Lambda$, it follows that necessarily
$K'=H\cap \Lambda$.
Hence $K$ contains $t_p$, and so $K=\langle \ell,t_p\rangle$.
This concludes the proof in the special case.

The general case is obtained by specialization.
Indeed, let $G_p\simeq \mathbb P^2$ be the space of lines through $p$. The subset
\[
\cW:=\{m\in G_p : m\neq t_p \text{ and } m\cap C=\{p\}\}
\]
is dense open in $G_p\smallsetminus\{t_p\}$, because its complement consists of the lines through $p$ meeting $C\smallsetminus\{p\}$, which form a proper closed subset. Since the family of limit cones $\widetilde C_0^m$ depends continuously on $m$, and the family of planes $\langle m,t_p\rangle$ clearly does as well, the identity
\[
\widetilde C_0^m=\widetilde C_p\cup \langle m,t_p\rangle
\]
proved for $m\in \cW$ extends to every $m\neq t_p$ by taking limits.
\end{proof}

\begin{remark}\label{rem: extension}
The same argument works in any dimension $n$ provided $C$ is an irreducible codimension 2 subvariety of $\bP^n$.
\end{remark}
\begin{remark}\label{rem: ciliberto}
The same result for a general line $\ell$ passing through $p$ can be derived by the arguments of Section 3 of \cite{CMR} (see in particular Lemma 3.1). 
In this part of the paper, the authors study more generally degenerations of projections of a variety $X\subset \pr^n$ when the centre of projection acquires a point of intersection with $X$.
\end{remark}
We also need to cover the case in which we take the limit along the tangent line. 
\begin{lemma}  \label{lem:osculating}
Let $p \in C$ and $t_p$ the projective tangent line at $p$. Then
$$\widetilde{C}^{t_p}_0= \widetilde C_p\cup H,$$ 
where $H$ is the osculating plane to $C$ at $p$. 
\end{lemma}
\begin{proof}
    Since the curve $C$ is smooth, A first-order expansion along the tangent direction shows that the residual plane contains $t_p$ and has contact order strictly greater than $2$ with $C$ at $p$. Hence it is the osculating plane to $C$ at $p$, see \cite{ACGH}.
\end{proof}

\subsection{Analytical proof of \Cref{lemma:limitconesexplicit}}\label{subs:coordinates}
We now give another proof of \Cref{lemma:limitconesexplicit}. This new proof involves an explicit computation, with a suitable choice of coordinates. Similar, but slightly more involved, computations are used in the arguments of \Cref{sec: conic series} and \ref{sec: conic divisors}.
Fix homogeneous coordinates $(x:y:z:w)$ of $\bP^3$. 
We can assume the coordinates of $p$ to be $(0:0:0:1)$, so that $W(p)=\langle x,y,z \rangle$, equivalently,  $f_p$ does not depend on the variable $w$.
We will also fix the tangent line to the curve $C$ at $p$ to be $t_p:=\{x=y=0\}$. Notice that with this choice, we have that  $\pi_p(p)=(0:0:1) \in \bP^2$.
Let us consider the line $\ell:=\{p_t=(-at:-bt:-ct:1), t \in \bC\}$ in $\bP^3$ as in \Cref{eqn:retta}. 
When moving the point $p_t$ along $\ell$, we get a deformation of $W(p)$ given by
$$W(p_t)=\langle x+atw,y+btw,z+ctw\rangle.$$ 
The  tangent vector to this deformation is
\begin{equation} \label{diff}
X=a\frac {\partial} {\partial x}+b\frac {\partial} {\partial y}+c\frac {\partial} {\partial z}.
\end{equation}
If we assume $\ell\neq t_p$, by taking into account all the previous assumptions, we may set $a=1, b=c=0$. Then, $p_t=(-t:0:0:1)$, $W(p_t)=\langle x+tw,y,z \rangle$. The line $\ell$ has equation $\{y=z=0\}$.

\begin{setting}\label{setting}
Here is an overview of all the choices we have discussed.
\begin{enumerate} \label{x}
\item $p=(0:0:0:1)$, a fixed point on $C$; 
\item $t_p=\{x=y=0\}$, the tangent line to $C$ at $p$; 
\item $\ell=\{y=z=0\}$ the implicit equations of the line along which we deform;
\item $p_t=(-t:0:0:1)$, the corresponding linear deformation of $p$; 
\item $W(p_t)=\langle x+tw,y,z\rangle$ the corresponding deformation of $W(p)$.
\end{enumerate}
\end{setting}
\begin{remark}\label{rem: nu}
Given any homogeneous polynomial $l\in V_k$ we will denote by $\nu_p(l)$ the evaluation of $l$ at $p$ on $C$. This is the intersection multiplicity of $C$ with the hypersurface $ Z(l)$ in the point $p$. 
 With the choices of \Cref{setting}, we have $\nu_p(w)=0$, $\nu_p(z)=1$, $\nu_p(x)\geq 2$ and $\nu_p(y)\geq 2$.
 \end{remark}

\begin{proof}[Proof of \Cref{lemma:limitconesexplicit}]
Without loss of generality, we assume Setting \ref{setting}.
Firstly we prove the result under the following assumptions, valid for a general $p \in C$ and $\ell$:
\begin{enumerate}
\item[($1^\star$)] $p \in C'$;
\item[($2^\star$)] $t_p \cap C = 2p \text{ scheme theoretically}$;
\item[($3^\star$)] $\ell$ is not bisecant to $C$.
\end{enumerate}
Notice that ($2^\star$) holds for the general point of $C$, see \cite{K} and \cite{BoP}.
From ($1^\star$) and \Cref{lem:suc} the projection curve $C_p=\pi_p(C)$ has degree $d-1$. From ($2^\star$), $\pi_p(p)$ is a smooth point of $C_p$.

Let $f_t=0$ be an equation of degree $d$ of the cone $\widetilde C_{p_t}$, for $t \neq 0$. The fact that $p_t$ is moving along $\{y=z=0\}$ can be encoded in the polynomial $f_{t}$ by highlighting the contribution of $t$ in the first variable: $f_t(x+tw,y,z)$. 

By \Cref{lemma: limitcones} the limit cone $\widetilde C_0^\ell$ has degree $d$. Since it contains the degree $d-1$ cone $\widetilde C_{p}$, it also contains a plane $H:=\{h=0\}$ passing through $p$. Since $h(p)=0$, and $p=(0:0:0:1)$, we have $h=Ax+By+Cz$. 
Let $f=0$ be an equation of  $\widetilde C_{p}$.
It follows that $f_0=h(x,y,z)f(x,y,z)$.

By assumption ($2^\star$) we have that the point $(1:0:0)$ does not belong to the tangent line of $C_p$ at $\pi_p(p)=(0:0:1)$. This is equivalent to 
 $\derx{f}|_{(0:0:1)} \neq 0$. 
 This means that, up to a non-zero constant, we can write 
\begin{equation} \label{eqn:nuovaf}
f=xz^{d-2}+ x^2L(x,y,z)+yM(x,y,z),
\end{equation}
where $\deg L=d-3$, $\deg M=d-2$.

The following formula follows from taking Taylor expansions in $t$, and it is proved in \Cref{lemma:taylor}. 
\begin{equation}\label{eq: taylor}
f_t(x+tw,y,z)-h(x,y,z)f(x,y,z)=t\left(w\derx{(h f)} +g(x,y,z)\right) \mod \ t^2.
\end{equation}
It follows that $w\derx{(h f)} =w(\derx{h} f+ h \derx{f})=wAf+ wh\derx{f}$, and so
\[
w\derx{(h f)} +g= wAf+ wh\derx{f}\in I(d).
\]
In particular we have
\begin{equation}\label{eq: derivata 2}
w h \derx{f}+g\in I(d).
\end{equation}
By evaluating \Cref{eq: derivata 2} at $p$, we obtain 
$$\nu_p(w)+\nu_p(h)+\nu_p\left( \derx{f}\right)=\nu_p(g)\geq d.$$ 
By \Cref{eqn:nuovaf}, we have  $\nu_p(\derx{f})=d-2$,
and so we obtain $\nu_p(h)\geq 2$. 
This means that 
$h(x,y,z)=Ax+By$, equivalently $H$ contains $t_p$. 
Moreover, $\ell \subset \widetilde C_{p_t}$ for all $t\not =0$, so $\ell \subset \widetilde C_p^\ell = \widetilde C_p\cup H$.

By Assumption ($3^\star$), we get $ \ell \subset H$, and then $H=\langle \ell,t_p \rangle$ is the plane spanned by $\ell$ and the tangent $t_p$, since $\ell=\{p_t=(-t:0:0:1), t\in \mathbb{C}\}$, we have that $ A=0 $ and $H= \{y=0\}$.
So,  \eqref{eq: derivata 2} becomes 
\begin{equation} \label{relation}
wy\derx{f}+g\in I(d),
\end{equation}

This proves the result under the assumptions ($1^\star$), ($2^\star$), ($3^\star$).
By continuity we can drop ($2^\star$) and ($3^\star$), and get the result for any $p\in C'$.

In order to drop ($1^\star$), consider a point $p \in S\cap C$. 
Choose a small neighbourhood  $\cV\subset C$ of $p$ such that any $q\in \cV\smallsetminus \{p\}$ is in $C'$. 
Consider any direction $v$ different from the one of $t_p$; then we can assume that $v$ is not the direction of  $t_q$ for any $q\in \cV$. Let $\ell$ be the line passing through $p$ with direction $v$. 
For any $q\in \cV\smallsetminus \{p\}$ we have that the  linear limit 
of $\widetilde{C}_{q_t}$ along a line with direction $v$ is the cone $\widetilde{C}_{q}\cup \langle t_q,\ell\rangle$ by the proof above. 
Then by specialization the same holds for $p$.
Note that in this case we have $\deg(\pi_p)=e$ and $er=d-1$. 
The equations of the limit cones have the form $f_p^e h$. 
As a divisor of $\bP^3$ we get $H+e\widetilde C_p$. This concludes the proof.
\end{proof}

\begin{lemma} \label{lemma:taylor}
With the assumptions and notations of \Cref{lemma:limitconesexplicit} we have
\begin{equation*}
f_t(x+tw,y,z)-h(x,y,z)f(x,y,z)=t\left(w\derx{(h f)} +g(x,y,z)\right) \mod \ t^2.
\end{equation*}
\end{lemma}
\begin{proof}
We highlight the dependence of $f_t(x+tw,y,z)$ on $t$ by writing $f_t=\varphi\circ \alpha$, where $\varphi\colon \bC^4\to \bC$ and $\alpha \colon \bC^5\to \bC^4$ is given by $\alpha (t,x,y,z,w)=(t, x+tw,y,z)$. 
Now we use the Taylor expansion of $\varphi\circ \alpha$  with respect to $t$ near $t=0$, calling $(x_1,\ldots , x_4)$ the coordinates in $\bC^4$:
\begin{align*}
\varphi\circ \alpha(t,x,y,z,w)&-\varphi\circ \alpha (0,x,y,z,w)=\\
&=t\frac{\partial \varphi}{\partial x_1} \alpha(0,x,y,z,w)+ tw\frac{\partial \varphi}{\partial x_2} \alpha(0,x,y,z,w) \mod \ t^2=\\
&=t\frac{\partial \varphi}{\partial x_1}(0, x,y,z)+ tw  \frac{\partial \varphi}{\partial x_2}(0, x,y,z) \mod \ t^2.
\end{align*}
Call $\frac{\partial \varphi}{\partial x_1}(0, x,y,z)=g(x,y,z)$, and observe that 
\begin{align*}
\frac{\partial \varphi}{\partial x_2}\alpha(0,x,y,z,w)&=\frac{\partial f_t}{\partial x}\alpha(0,x,y,z,w)\frac{\partial \alpha}{\partial x_2}(0,x,y,z,w)=\\
&=\frac{\partial f_0}{\partial x}(x,y,z)=\frac{\partial}{ \partial x}((h\cdot f)(x,y,z)), 
\end{align*}
since $\frac{\partial \alpha}{\partial x_2}(0,x,y,z,w)= 1$, formula \eqref{eq: taylor} is verified.
\end{proof}

\begin{remark}\label{limit S}
Consider the case where $p\in S$, $p\notin C$. 
Assume that $\deg \pi_p =e$ and $d=er$.
Let $\widetilde C_p$ be the cone of degree $r$ joining $p$ and $C$. 
If we take any limit $p_t\to p$, the equation of the limit cone becomes
$f_p^e$, where $f_p\in I(r)$ is the equation of the reduced cone, and the divisor we get in $\bP^3$ is then $\widetilde C_p$ with multiplicity $e$. 
\end{remark}

\subsection{A projective model of the blow-up of $\bP^3$ along $C$}\label{sec: proj model}
Let $\nu\colon \widetilde \bP^3\to \bP^3$ be  the blow-up of $\bP^3$ along $C$. We  now extend the map (\ref{c map}) over $\widetilde \bP^3$.
\begin{proposition}\label{extension c}
The map $c\colon U\to \bP (I(d))$  extends to  an injective morphism $$\widetilde c\colon \widetilde \bP^3\longrightarrow \bP (I(d)).$$
The closure of $c(U)\subset \bP (I(d))$ is $\widetilde c(\widetilde \bP^3)$.
The divisor associated to $\widetilde c$ in $\pic(\widetilde \bP^3)$ is $d\widetilde L-E$, where $\widetilde L$ is the pullback of a hyperplane divisor on $\pr^3$ and $E$ is the exceptional divisor of $\nu$.
\end{proposition}
\begin{proof}
The exceptional divisor of $\nu$ is precisely $E=\pr(N_{C|\bP^3})$.
The extension of $c$ can be derived from Lemma \ref{lemma:limitconesexplicit} as follows. 
Recall that any $\xi\in E$  corresponds to a point $p\in C$ and to an element in the fibre of the projectivized of the normal bundle $\pr(N_{C|\bP^3})_{p'}\cong \pr^1$, which corresponds to a plane $H$ passing through $p$ and containing $t_p$. 
Define $\widetilde{c}(\xi):=[f_{p}h]\in \bP (I(d))$, where $\{h=0\}$ is the equation of $H$. 
On the finite set $S\smallsetminus C$ the limit is unique as observed in Remark \ref{limit S}, so the extension is well defined on the whole $\widetilde \bP^3$. 
It is easy to prove that it is injective. 
Let us just observe that in the case we have two distinct points $\xi,\xi'\in E$ such that $p=p'$, then the hyperplanes are different, and so $\tilde c(\xi)\not=\tilde c(\xi')$.

We now want to compute the class of ${\widetilde c}^{*}\cO_{\bP (I(d))}(1)$. Recall that 
$$\pic(\widetilde \bP^3)\cong \bZ[\widetilde L]\oplus\bZ[E],$$ 
where $\widetilde L$ is the pullback of a hyperplane divisor $L$ of $\bP^3$ via $\nu$.

Following for instance the results of \cite[Chap.4, sec. 6]{GH}, it can be easily proved that:
\[
\widetilde L^3= L^3=1, \quad
\widetilde L^2\cdot E=0, \quad
\widetilde L\cdot E^2=-d,\quad\quad\quad\quad\quad\quad\]
\[E^3=-\deg N_{C|\pr^3}=-2g+2+(K_{\pr^3}\cdot C)=-2g+2-4d.
\]
Observe that an effective divisor in $|{\widetilde c}^{*}\cO_{\bP (I(d))}(1)|$ can be seen as ${\widetilde c}^{\ast}\Sigma$, where: 
\[
\Sigma:=\{\mbox{classes of hypersurfaces of degree } d \mbox{ containing } r\mbox { and } C\}\subset \bP(I(d)),
\]
where $r\in \bP^3$ is a general point.

Let $\alpha, \beta \in \bZ$ such that $|{\widetilde c}^{*}\cO_{\bP (I(d))}(1)|= |\alpha \widetilde L-\beta E|$. 
Call $M$ a general ${\widetilde c}^{\ast}\Sigma\in |{\widetilde c}^{*}\cO_{\bP (I(d))}(1)|$.
We prove the following two facts: 
\smallskip

\noindent  (1) $M\cdot \widetilde L\cdot E=d$.\\
Indeed, with the notation introduced above, the elements in $M\cap E$ are the pullbacks via $\widetilde{c}$ of the degenerate cones $\widetilde{C}_p\cup H$ where $p\in  C$, $p \in H$ and $r\in  \widetilde C_p\cup H$. 
We now intersect with $\widetilde L=\nu^\ast L$. The intersection $L\cap C$ is a set of $d$ points in general position by the general position theorem \cite[Chap.2 Sec.3]{GH}, which we call $p_1,\ldots p_d$. 
For a general choice of $L$ and $r$ we have that $r$ does not belong to any of the cones $\widetilde C_{p_i}$ with vertex $p_i$ for any $i=1,\ldots, d$. 
Therefore, the elements in $M\cdot \widetilde L\cdot E$ are the pullback of the set of degenerate cones $\widetilde C_{p_i}\cup H_i$ such that  $r\in H_i$. 
Recall that $H=\langle t_{p_i}, \ell\rangle $ for some line $\ell\subset \bP^3$. 
But for any $i$ there exists exactly one line $\ell_i$ such that $r\in \langle t_{p_i}, \ell_i\rangle $. 
So, for any $i\in \{1,\ldots ,d\}$ we have precisely one contribution to  $M\cdot \widetilde L\cdot E$, and the proof is concluded.

\smallskip

\noindent (2) $M^2\cdot E=2\left((d-1)^2-g\right)$.\\
 Indeed, $M^2$ can be seen as the pullback via $\widetilde c$ of the hypersurfaces in $I(d)$ containing two general points $r, s$ in $\bP^3$. 
 Since $r, s$ are general, then they never are contained in the same cone $\widetilde C_p$ for $p\in \bP^3\smallsetminus C$. 
 On the other hand, there are $2g-2+2d$ degenerate cones $\widetilde{C}_p\cup H$ (whose pullback is in $E$) such that the plane $H$ contains both $r$ and $s$. Indeed, consider the projection from the line $\langle r,s\rangle$:
\[\xymatrix{
\bP^2\smallsetminus \langle r,s\rangle\ar[r]& \bP^1\\
C\ar@{^{(}->}[u]\ar[ru]&
}\]  
Its restriction to $C$ is a degree $d$ covering with only simple ramifications by the generality of $r$ and $s$. 
The cones we want to count correspond to these points of ramifications, so we obtain $2g-2+2d$ such points by the Riemann-Hurwitz formula \cite[Preliminaries]{GH}.

Now we are left to compute the number of degenerate cones $\widetilde{C}_p\cup H$ such that $H$ contains $r$ and $\widetilde{C}_p$ contains $s$.

Observe that $s\in \widetilde{C}_p$ if and only if the line $\langle p,s\rangle$ intersects the curve $C$. In other words, the image curve obtained by  projecting from $s$, $\pi_s\colon C\to C_s\subset\bP^2$ has a node.  
Now, the number of nodes of $C_s$ coincides with the difference between its arithmetic and its geometric genus:
\[
p_a(C_s)-g=\frac{(d-1)(d-2)}{2}-g.
\]
So, for every node $q$ of $C_s$ we have two points in $ M^2\cap  E$: let $\{p,p'\}=\pi_s^{-1}(q)$, we have then precisely the two cones 
\[\widetilde C_p\cup \langle t_p, \ell_{s,p}\rangle, \quad \widetilde C_{p'}\cup \langle t_{p'}, \ell_{s,p'}\rangle.
\] 
Now, switching the roles of $r$ and $s$ we have another equal contribution.
Summing up, we have 
\[
M^2\cdot E= (2g-2+2d) + 2\big((d-1)(d-2)-2g\big)=2((d-1)^2-g),
\]
as wanted. 

Now, let us use the information above. We have from (1) and the relations listed above  that 
\[
-d=M\cdot \widetilde L\cdot E=(\alpha \widetilde H-\beta E)\cdot \widetilde L\cdot E=-d\beta,
\]
so $\beta =1$. 
From (2) and what we have  proven so far we have:
\[
2((d-1)^2-g)=M^2\cdot E=(\alpha \widetilde L- E)^2\cdot E=2\alpha d+2-2g-4d,
\]
so we obtain $\alpha=d$ as wanted.
\end{proof}

\begin{remark}
So we have by \cite[Cor.1.2.15]{Laz1} that $d\widetilde L-E$ is ample. 
Note that $-K_{\widetilde{\bP^3}}=4\widetilde L-E$. 
In case $d=3,4$ it is already known that $-K_{\widetilde{\bP^3}}$ is ample \cite[Prop. 3.1]{BL}.  
\end{remark}

\section{Conic linear systems and the cone map}\label{sec: conic series}
Now we want to consider divisors on $C$ coming from cones in $\bP^3$. For any integer $k\geq 1$, call $S_k:=\H^0(C,\cO_C(k))$, so that $\bP (S_k)$  parametrises divisors of degree $dk$ on $C$ cut out on $C$ by hypersurfaces of degree $k$ in $\pr^3$.
Recall that $V_k=H^0(\pr^3,\cO_{\pr^3}(k))$; denote by $\phi_k \colon V_k\to S_k$ the natural restriction map. 
We have the sequence 
\[
0\longrightarrow I(k)\longrightarrow V_k\stackrel{\phi_k}{\longrightarrow} S_k,
\]
where $I(k)=\ker \phi_k$ is the space of the homogeneous polynomials of degree $k$ vanishing on $C$. 
We know from the Castelnuovo type result in \cite{GLP} that $\phi_k$ is surjective for $k\geq d-2$ and from Castelnuovo's results \cite{Cast}, \cite[First Theorem of Castelnuovo]{szip} that $H^1(C,\cO_C(k))=0$ for $k\geq \lfloor\frac{d-1}{2}\rfloor$. 

Let us assumet that $k\geq d-2$.
By applying the Riemann-Roch theorem, we get 
$\dim S_k = dk-g+1$, and in particular
$\dim S_d= d^2-(g-1)$.
It follows that $I(d)=\ker \phi_d$ has dimension $\binom{d+3}{3}-d^2+g-1$.  

\begin{definition} \label{limnormali} 
Given a linear subspace $W\subset V$, the {\em conic linear system} over $C$ of degree $dk$ associated to $W$ is 
$$R_k(W) := \phi_k (W_k)\subset S_k.$$ 
For $p\in \bP^3$ we will use the notation $R_k(p)=R_k(W(p))=\phi_k (W(p)_k)$.
 \end{definition}

We also call \emph{conic linear series} the projectivization of a conic linear system and \emph{conic divisor} an element of a conic linear series.
 
 \begin{remark}\label{rem: dimR}
Consider the following diagram: 
\[\xymatrix{
&0\ar[d]&0\ar[d]&0\ar[d]&\\
0\ar[r]& I(k)\cap W_k\ar[d]\ar[r]& W_k\ar[d]\ar[r] &R_k\ar[d]\ar[r]& 0\\
0\ar[r]& I(k)\ar[r] &V_k\ar[r] &S_k&\\
}\]
Let us consider $W=W(p)$. For $k=d-1$ we have:
\begin{itemize}
\item[(1)] $p\in  U$ if and only if $\dim R_{d-1}(p)=\dim W(p)_{d-1}=\binom{d+1}{2}$;
\item[(2)] $p\in  C'$ if and only if $\dim R_{d-1}(p)=\dim W(p)_{d-1}-1=\binom{d+1}{2}-1=\frac{(d-1)(d+2)}{2}$;
\item[(3)] $p\in  S$ if and only if  $\dim R_{d-1}(p)\leq \dim W(p)_{d-1}-2=\binom{d+1}{2}-2=\frac{d^2+d-4}{2}$.
\end{itemize}
By \Cref{lem: dimensioni}, for $k=d$ we have instead:
\begin{itemize}
\item[(1)] $p\in  U$ if and only if $\dim R_d(p)=\dim W(p)_d-1=\binom{d+2}{2}-1=\frac{d(d+3)}{2}$;
\item[(2)] $p\in  C'$ if and only if $\dim R_d(p)=\dim W(p)_d-3=\binom{d+2}{2}-3=\frac{(d-1)(d+4)}{2}$;
\item[(3)] $p\in  S$ if and only if $\dim R_d(p)\leq\dim W(p)_d-6=\binom{d+2}{2}-6=\frac{d^2+3d-10}{2}$.
\end{itemize}

\end{remark}

For a vector space $A$ and an integer $s\leq \dim A$ we let $G(s,A)$ be the Grassmannian of the $s$-vector spaces of $A$. 
\begin{definition}\label{def: rho}
Let $k\leq d$ and take $p \in U$. We have a morphism
\begin{align}\label{rho}
\begin{split}
    \rho_k\colon U&\longrightarrow G(n_k,S_k)\\
    p&\longmapsto R_k(p)\subset S_k
    \end{split}
\end{align}
where $n_k:=\dim R_k(p)$.
We call this map {\em the $k$-th cone map}.
\end{definition}
By the above discussions and by \Cref{rem: dimR} we have that $n_k=\binom{k+2}{2}$ for $k<d$, and that $n_d=\binom{d+2}{2}-1$.
\begin{remark}\label{rem: cone map}
The cone map associates to any $p\in U$ the conic linear system which is the image  $R_k(p)$ of $W(p)_k$ in $S_k$. 
Observe 
that in case $k=d$, we have that $R_d(p)$ is isomorphic to $W(p)_d/\langle f_p\rangle$, where $f_p$ is as usual an equation (unique up to $\mathbb C^*$) of the cone over $C$ with vertex in $p$.
\end{remark}
\begin{proposition}\label{lem: inj rho}
The morphism $ \rho_d$ is generically injective for $d\geq 4$. 
\end{proposition}
\begin{proof}
Assume that $p$ and $q$ are distinct points of $U$ and $\rho_d(p)=\rho_d(q)$. 
As a consequence we have that $R_d(p)=R_d(q)$, i.e. the linear systems cut out by $W(p)_d$ and $W(q)_d$ on $C$ are the same. 
In particular, the points that are identified by the morphisms $\pi_p\colon C\to C_p\subset \bP(W(p)^\ast)$ and $\pi_q\colon C\to C_q\subset \bP(W(q)^\ast)$ are the same.
Now, two points of $C$ are identified by  $\pi_p$ if and only if  there exists a line $r$ passing through $p$ intersecting $C$ in two distinct points. If $\rho_d(p)=\rho_d(q)$ necessarily $\pi_q$ identifies these points, so $q$ belongs to $r$.
For $d\geq 4$ and $p $ general, there exist at least two bisecants $r$, $s$ passing through $p$, because the image curve is nodal and has at least two nodes.  So we have that $q\in r\cap s=\{p\}$. 
\end{proof}
We have just proved the general injectivity of $\rho_d$ for $d\geq 4$ by using that fact that there exist at least two bisecants passing from a general point of a curve in $\bP^3$ of degree at least $4$. 
This last statement is false for the rational normal cubic in $\bP^3$. Indeed, we have the following result:
\begin{lemma}
Let $C\subset \pr^3$ be the rational normal cubic. The morphism $\rho_3$ is generically $3:1$. 
\end{lemma}
\begin{proof}
Fix a general point $p \in U$. From $p$ there exists a unique bisecant to $C$, call it $r$. We have seen in the proof of \Cref{lem: inj rho} that any point $q\in U$ such that $\rho_3(p)=\rho_3(q)$ has to lie on $r$. Now we see that there are three points in the preimage $\rho_3^{-1}(\rho_3(p))$. Let us fix the coordinates so that $p=(1:0:0:-1)$, and consider $p':=(a:0:0:-1)$, $p''=(a^2:0:0:-1)$, where $a$ is a primitive cubic root of $1$. 
We have that 
$\rho_3(p)=\rho_3(p')=\rho_3(p'')$. 
Indeed, observe first that 
\[
W(p)=\langle y,z, x+w\rangle,\,\, W(p')=\langle y,z, x+aw\rangle, \,\,W(p'')=\langle y,z, x+a^2w\rangle.
\]
Let us consider the affine open set $w=1$; if we substitute $x=t^3, y=t^2, z=t, w=1$ in the generators of $W(p)_3$, $W(p')_3$ and $W(p'')_3$ respectively, we obtain in any case 
\[R_3(p)=R_3(p')=R_3(p'')=\langle
t,t^2,t^3,\ldots,t^8,t^9+1
\rangle.
\]
If we consider $q=(b:0:0:-1)$, with $b^3\not =1$, it is easy to see that over $C$ we obtain 
\[
R_3(q)=\langle t,t^2,t^3,\ldots,t^8,t^9+b^3 \rangle\not= R_3(p).
\]
So the proof is concluded. 
\end{proof}

\section{Limits of conic linear systems and extension of the cone map}\label{sec: limits}
As above, let $\nu:\widetilde{\mathbb P^3}\longrightarrow \mathbb P^3$
be the blow-up of $\mathbb P^3$ along $C$, and let $E=\mathbb P(N_{C|\mathbb P^3})$
be its exceptional divisor. 
A point of the exceptional divisor over $p\in C$ will be denoted by $\xi=(p,H)\in E$, where $H$ is a plane containing the tangent line $t_p$. 
We also call $\mu\colon E\longrightarrow C$ the natural projection.

We have the cone map
$$
\rho_k\colon U\longrightarrow G(n_k,S_k).
$$
In this section we extend $\rho_k$ to a rational map $\widetilde{\rho}_k\colon \widetilde{\mathbb P^3}\dashrightarrow G(n_k,S_k),$
and we give geometric information to the maximal open subset of $\widetilde \pr^3$ where it can be defined. 
The naif idea would be to proceed as in \Cref{sec: proj model}, just considering linear limits of conic series, in order to extend $\rho_k$ to the exceptional divisor $E$. 
But -as it was pointed out to us by Wei Chen- it is not clear that this limits only depends on a class $\xi\in E$ for $k=d-1,$ or $k=d$.
So, we will have to make a more delicate construction: first we extend formally $\rho_k$ to a codimension at least 2 open set, then we use analytic computations to understand the intersection of this open set with $E$.

By composing with the Pl\"ucker embedding of $G(n_k,S_k)$, we may regard ${\rho}_k$ as a rational map to a projective space. 
We can thus extend $\rho_k$ to  a map
$$\widetilde{\rho}_k\colon \widetilde{\mathbb P^3}\dashrightarrow G(n_k,S_k)$$
whose indeterminacy locus has codimension at least $2$. 
We denote by $\mathcal V_k\subseteq \widetilde{\mathbb P^3}$ the complement of the indeterminacy locus, i.e.
the maximal open subset on which $\widetilde{\rho}_k$ is regular. 
Note that $U \subset \mathcal V_k$.

\begin{definition}\label{def: Rxi-regular}
For an element $\xi \in \mathcal V_k\cap E$, with $p=\mu(\xi)$, we define
$$
R_k^\xi(p):=\widetilde{\rho}_k(\xi).
$$
\end{definition}

We now relate this abstract definition with limits computed along lines. 
Let $p\in C$, and let $\ell\neq t_p$ be a line through $p$. 
Assume that $\ell\cap U\neq \emptyset$. 
Then the restriction
$$
{\rho_k}_{|\ell\cap U}\colon \ell\cap U\longrightarrow G(n_k,S_k)
$$
extends uniquely to $\ell\simeq \mathbb P^1$, since the Grassmannian is projective. We denote by $R_k^\ell(p)$ the value of this extension at $p$.

\begin{lemma}\label{lem:line-computes-rho}
Let $\xi=(p,H)\in \mathcal V_k\cap E$, and let $\ell$ be a line through $p$, with $\ell\neq t_p$ and $\langle t_p,\ell\rangle=H$. Assume that $\ell\cap U\neq \emptyset$. Then $R_k^\ell(p)=R_k^\xi(p).$
\end{lemma}
\begin{proof}
Let $p_t\in \ell\cap U$ be a point that tends to $p$, and let $\widetilde p_t$ be its lift to $\widetilde{\mathbb P^3}$. Since $\ell$ represents the normal direction $\xi=(p,H)$, we have $\widetilde p_t\to \xi$. Hence, by the regularity of $\widetilde\rho_k$ at $\xi$,
\[
R_k^\ell(p)=\lim_{t\to 0}R_k(p_t)
=\lim_{t\to 0}\widetilde\rho_k(\widetilde p_t)
=\widetilde\rho_k(\xi)=R_k^\xi(p).
\]
\end{proof}

\begin{lemma}\label{lem: inclusione}
Consider a point $p \in C$ contained in a line $\ell \neq t_p$ such that $\ell \cap U \neq \emptyset$. There is an inclusion $R_k(p)\subseteq R^\ell_k(p).$
\end{lemma}
\begin{proof}
Choose coordinates as in \Cref{setting}. Fix a cone of degree $k$ given by a polynomial $r(x,y,z)\in W(p)_k$.
If we consider
$r(x+tw,y,z)\in W(p_t)_k$,
for $t\in \mathbb C$, we have a family of cones with vertex $p_t$ whose limit is $r(x,y,z)$. Therefore
$$
\phi_k(W(p)_k)\subseteq \lim_{t\to 0}\phi_k(W(p_t)_k),
$$
as wanted.
\end{proof}

\begin{remark}\label{rmk:boh d-1}
In particular, if $k<d-1$ and $p\in C'$, we have
$R_k(p)=R^\ell_k(p)$,
because they have the same dimension.
\end{remark}

We now study the behaviour of $\widetilde{\rho}_k$ on the exceptional divisor for $k=d-1$ and $k=d$. 
Let
\[
B_k:=E\smallsetminus (E\cap \mathcal V_k)
\]
be the locus in $E$ where the value $R_k^\xi(p)$ is not defined by the regularity of $\widetilde{\rho}_k$. Since the complement of $\mathcal V_k$ has codimension at least $2$ in $\widetilde{\mathbb P^3}$, the locus $B_k$ has dimension at most $1$. Therefore, a priori, $B_k$ could contain vertical components, namely fibres of $\mu\colon E\to C$, and horizontal components dominating $C$.

We will prove that the possible indeterminacy of $\widetilde\rho_{d-1}$ and $\widetilde\rho_d$ along $E$ is concentrated on fibres over a precise finite subset of $C$. Let
\[
\Gamma:=(C\smallsetminus C')\cup 
\{p\in C' \mid t_p\cap C\neq \{p\}\}.
\]
Equivalently, $\Gamma$ is the finite set of points of $C$ where at least one of the following two assumptions fails:
\begin{enumerate}
\item[($1^\star$)] $p\in C'$;
\item[($2^\star$)] $t_p \cap C = 2p \text{ scheme theoretically}$.
\end{enumerate}
The finiteness of $\Gamma$ follows from \Cref{lem:suc} and from the results recalled in the proof of \Cref{lemma:limitconesexplicit}. 

In the resto  of the section we will prove the following result.
\begin{theorem}\label{teo: VsuE}
With the above notations, for  $k=d-1$ and $k=d$, one has
\[
\mu(B_k) \subset \Gamma.
\]
Equivalently, for every $p\in C\smallsetminus\Gamma$, the system $R_k^\xi(p)$ is defined for every normal direction $\xi\in \mu^{-1}(p)$.
\end{theorem}

\subsection*{The case of $k=d-1$}

Fix $p\in C\smallsetminus \Gamma$ and $\xi=(p,H)\in E$. We choose coordinates as in \Cref{setting}. The lines through $p$ contained in $H$ and different from $t_p$ are written as
\[
\ell_a:=\{y=0,\ z-ax=0\}.
\]

\begin{lemma}\label{lem:line-limit-d-1}
Let $p\in C\smallsetminus\Gamma$, and let $\xi=(p,H)\in E$. With the notation above, for every line $\ell_a\subset H$ we have
\[
R^{\ell_a}_{d-1}(p) = \left\langle
R_{d-1}(p),
w\left(f_x-af_z\right)
\right\rangle .
\]
In particular, we have $\nu_p(w(f_x-af_z))=d-2$.
\end{lemma}
\begin{proof}
We deform the vertex along $\ell_a$ by setting
$p_t=(-t:0:at:1)$. Then $W(p_t)=\langle x+tw,\ y,\ z-atw\rangle$.
By \Cref{lem: inclusione}, we already have
$R_{d-1}(p)\subseteq R^{\ell_a}_{d-1}(p)$.

Consider the first order expansion
\[
f(x+tw,y,z-atw)
=
f(x,y,z)+t\,w(f_x-af_z)\mod t^2.
\]
After restricting to $C$, the term $f(x,y,z)$ vanishes. Hence the first non-zero term gives
\[
w(f_x-af_z)\in R^{\ell_a}_{d-1}(p).
\]
Therefore
\[
\left\langle
R_{d-1}(p),\,
w(f_x-af_z)
\right\rangle
\subseteq
R^{\ell_a}_{d-1}(p).
\]

It remains to show that the inclusion is an equality. Since $p\in C\smallsetminus\Gamma$, the assumptions $(1^\star)$ and $(2^\star)$ hold, and so $f$ can be written in the form of \eqref{eqn:nuovaf}:
\[
f=xz^{d-2}+x^2L(x,y,z)+yM(x,y,z),
\]
with $\deg L=d-3$ and $\deg M=d-2$. Thus
\begin{align*}
    f_x&=z^{d-2}+\text{terms of evaluation at least }d-1 \text{ in }p,\\
    f_z&=(d-2)xz^{d-3}+\text{terms of evaluation at least }d-1 \text{ in }p.
\end{align*}
Since $\nu_p(w)=0$, $\nu_p(z)=1$ and $\nu_p(x)\geq 2$, we get
$\nu_p(wf_x)=d-2$, $\nu_p(wf_z)\geq d-1$.
Therefore, for every $a$,
\[
\nu_p\big(w(f_x-af_z)\big)=d-2.
\]
On the other hand, every element of $R_{d-1}(p)$ has evaluation at least $d-1$ at $p$. and so $w(f_x-af_z)\notin R_{d-1}(p)$.
Thus
\[
\dim \left\langle
R_{d-1}(p),\,
w(f_x-af_z)
\right\rangle
=
\dim R_{d-1}(p)+1.
\]
Since $R^{\ell_a}_{d-1}(p)$ is a limit in $G(n_{d-1},S_{d-1})$, it has dimension $n_{d-1}$. Moreover, as $p\in C'$, by \Cref{rem: dimR} we have $\dim R_{d-1}(p)=n_{d-1}-1$.
The inclusion above is therefore an equality:
\[
R^{\ell_a}_{d-1}(p)
=
\left\langle
R_{d-1}(p),\,
w(f_x-af_z)
\right\rangle .
\]
\end{proof}

\begin{lemma}\label{lem:L-equals-R}
Let $\xi=(p,H)\in \mathcal V_{d-1}\cap E$, with $p\in C\smallsetminus \Gamma$. Then
\[
R_{d-1}^\xi(p) = \left\langle
R_{d-1}(p),
wf_x
\right\rangle .
\]
In particular, this description does not depend on the line $\ell$ chosen to represent $\xi$.
\end{lemma}
\begin{proof}
Since $\xi\in \mathcal V_{d-1}$, by \Cref{lem:line-computes-rho} we have $R_{d-1}^\xi(p)=R_{d-1}^{\ell_a}(p)$
for every line $\ell_a\subset H$ representing $\xi$. Hence, by \Cref{lem:line-limit-d-1},
$$
R_{d-1}^\xi(p)
=\left\langle
R_{d-1}(p),
w(f_x-af_z)
\right\rangle .
$$
Let
$$L^\xi:=\{u\in R_{d-1}^\xi(p)\mid \nu_p(u)\geq d-1\}.$$
The subspace $R_{d-1}(p)$ is contained in $L^\xi$. 
By the evaluation described in \Cref{lem:line-limit-d-1}, $L^\xi$ is a proper subspace of $R_{d-1}^\xi(p)$. Since $\dim R_{d-1}^\xi(p)=\dim R_{d-1}(p) + 1$, we get
$$
L^\xi=R_{d-1}(p).
$$
Now take two lines $\ell_a,\ell_b\subset H$. 
Since both the extra generators belong to $R_{d-1}^\xi(p)$, also $w(f_x-af_z)-w(f_x-bf_z)=(b-a)wf_z$ does, and so we have $wf_z\in R_{d-1}^\xi(p)$. Moreover $\nu_p(wf_z)\geq d-1$.
Therefore $wf_z\in L^\xi=R_{d-1}(p)$. Consequently
$$
w(f_x-af_z)\equiv wf_x\mod R_{d-1}(p),
$$
and so
$$
R_{d-1}^\xi(p)=\left\langle
R_{d-1}(p), wf_x
\right\rangle .
$$
\end{proof}

We now use this computation to extend the definition also at the points where the regularity of $\widetilde{\rho}_{d-1}$ is not known a priori.

\begin{proposition}\label{prop:no-horizontal-d-1}
With $\Gamma$ as above, we have
\[
\mu(B_{d-1}) \subset \Gamma.
\]
Equivalently, for $p\in C\smallsetminus \Gamma$ and every $\xi=(p,H)\in E$, the limit system $R_{d-1}^\xi(p)$ is well-defined.
\end{proposition}

\begin{proof}
Suppose, by contradiction, that there exists a point
\[
\xi=(p,H)\in B_{d-1}
\]
with $p\in C\smallsetminus\Gamma$. Since $p\notin\Gamma$, the assumptions $(1^\star)$ and $(2^\star)$ hold, and we may use the local form \eqref{eqn:nuovaf}.
By \Cref{lem:line-limit-d-1}, for every $a$ we have
\[
R_{d-1}^{\ell_a}(p)
=
\left\langle
R_{d-1}(p),w(f_x-af_z)
\right\rangle .
\]

We now prove that $wf_z\in R_{d-1}(p)$. Since $B_{d-1}$ has dimension at most $1$ in the surface $E$, we can choose a curve
\[
\xi_t=(p_t,H_t)\in (\mathcal V_{d-1}\cap E)
\]
such that $\xi_t\to \xi$, with $p_t\in C\smallsetminus\Gamma$ for $t\neq 0$. For every $t\neq 0$, by \Cref{lem:L-equals-R}, the subspace of elements of $R_{d-1}^{\xi_t}(p_t)$ with evaluation at least $d-1$ is precisely $R_{d-1}(p_t)$.

Choose coordinates adapted to $\xi_t=(p_t,H_t)$, and write $f_{z,t}$ for the derivative of the equation of $\widetilde C_{p_t}$ with respect to the corresponding coordinate $z_t$. By the same argument used in the proof of \Cref{lem:L-equals-R}, we have $wf_{z,t}\in R_{d-1}(p_t)$. As $t\to 0$, we have $wf_{z,t}\longrightarrow wf_z$.
Moreover, since $p_t\to p$ inside $C'$, the spaces $R_{d-1}(p_t)$ tend to $R_{d-1}(p)$ in the Grassmannian. Hence $wf_z\in R_{d-1}(p)$.

Therefore
\[
w(f_x-af_z)\equiv wf_x \mod R_{d-1}(p)
\]
for every $a$. Consequently
\[
R_{d-1}^{\ell_a}(p)
=
\left\langle
R_{d-1}(p),wf_x
\right\rangle
\]
for every line $\ell_a\subset H$ through $p$, different from $t_p$.

Thus the limit along a line contained in $H$ is independent of the chosen line. We may therefore define
\[
R_{d-1}^\xi(p):=
R_{d-1}^{\ell}(p)
=
\left\langle
R_{d-1}(p),wf_x
\right\rangle ,
\]
where $\ell\subset H$ is any line through $p$ different from $t_p$.

The computation above applies not only to linear deformations. Indeed, in adapted local coordinates, any curve in $\bP^3$ whose strict transform passes through $\xi=(p,H)$ has the same first-order expansion as one of the lines $\ell_a\subset H$, up to terms of order at least $2$. These higher-order terms do not affect the first non-zero contribution modulo $R_{d-1}(p)$. Hence every curve approaching $\xi$ gives the same limit
$$
\left\langle R_{d-1}(p),wf_x\right\rangle.
$$
Moreover, this subspace varies regularly with $(p,H)$ in a local chart of the blow-up. Therefore it defines a regular extension of $\widetilde\rho_{d-1}$ at $\xi$.

The construction above gives an extension of $\widetilde\rho_{d-1}$ at $\xi$. This contradicts the maximality of the open subset $\mathcal V_{d-1}$.
\end{proof}

\begin{corollary}\label{cor:evaluation-d-1}
Let $p\in C\smallsetminus \Gamma$ and let $\xi=(p,H)\in E$. Then every element
$u\in R_{d-1}^\xi(p)$ 
satisfies 
$\nu_p(u)\geq d-2.$
\end{corollary}
\begin{proof}
By the explicit description above,
$$
R_{d-1}^\xi(p)
=\left\langle
R_{d-1}(p),wf_x
\right\rangle .
$$
Every element of $R_{d-1}(p)$ has evaluation at least $d-1$, while
$\nu_p(wf_x)=d-2$.
The claim follows.
\end{proof}

\subsection*{The case $k=d$}

Fix $p\in C\smallsetminus\Gamma$ and $\xi=(p,H)\in E$. We make the same
coordinate choices as in the previous case. Thus $H=\{y=0\}$, the lines through
$p$ contained in $H$ and different from $t_p$ are
$\ell_a=\{y=0,\ z-ax=0\}$, and the corresponding deformation is
$$
W(p_t)=\langle x+tw,\ y,\ z-atw\rangle.
$$

Since $\dim R_d(p)=n_d-2$, two additional elements occur in the limit. 
Unlike the case $k=d-1$, the two additional elements may appear at different orders in $t$, making their explicit description more involved. Since in what follows we only need their evaluations at $p$, we will not describe them explicitly, but rather compute their evaluations from the first non-zero terms of the corresponding expansions after restriction to $C$.

\begin{lemma}\label{lem:two-step-limit-d}
There are classes
$\sigma_1^\xi,\sigma_2^\xi\in S_d$, varying regularly with $\xi$ locally on
$E$, such that, for every line $\ell_a\subset H$ representing $\xi$,
$$
R_d^{\ell_a}(p)
=
\left\langle
R_d(p),\sigma_1^\xi,\sigma_2^\xi
\right\rangle.
$$
Moreover,
$
\nu_p(\sigma_1^\xi)=d-1$, and $\nu_p(\sigma_2^\xi)=d-2$.
\end{lemma}

\begin{proof}
Let $f\in W(p)_{d-1}$ be an equation of the cone $\widetilde C_p$. Then
$$
K_p:=W(p)_d\cap I(d)=\langle xf,yf,zf\rangle.
$$
For $Q\in K_p$, write
$$
Q(x+tw,y,z-atw)
=
Q+t\,w(Q_x-aQ_z)\mod t^2.
$$
After restriction to $C$, the first-order coefficient defines a linear map
\begin{align*}
\delta_a: K_p&\longrightarrow S_d/R_d(p)\\
Q&\longmapsto
\left[\phi_d\bigl(w(Q_x-aQ_z)\bigr)\right]
\end{align*}
A direct computation on the basis $xf,yf,zf$ gives
$\operatorname{rank}\delta_a=1$.
Its image is generated by the class obtained from $zf$, $wz(f_x-af_z)$, and the local form \eqref{eqn:nuovaf} gives evaluation $d-1$.

Let $\mathcal L_a\subset\ker\delta_a$ be the line obtained as the limit of the
one-dimensional spaces $W(p_t)_d\cap I(d)$ for $t\neq0$. If
$Q\in\ker\delta_a$, choose $G\in W(p)_d$ cancelling the first-order term.
The second-order coefficient, considered modulo
$\left\langle R_d(p),\operatorname{Im}\delta_a\right\rangle$,
is independent of the choice of $G$ and depends only on the class of $Q$
modulo $\mathcal L_a$. It therefore defines a linear map
$$
\epsilon_a\colon
\ker\delta_a/\mathcal L_a
\longrightarrow
\frac{S_d}
{\left\langle R_d(p),\operatorname{Im}\delta_a\right\rangle}.
$$
The source is one-dimensional, and again a direct computation gives $\operatorname{rank}\epsilon_a=1$,
with non-zero image of evaluation $d-2$.

The constructions of $\delta_a$ and $\epsilon_a$ are the first two graded pieces of the saturated module of restrictions on a local chart of the blow-up. 
Since both have constant rank $1$, the first construction adds one elemetn to $R_d(p)$ and the second adds one further independent element. 
These elements vary regularly with $\xi$, so together with $R_d(p)$ they define a subbundle of rank $n_d$ and hence a local extension of $\widetilde\rho_d$.
\end{proof}

\begin{remark}
The classes
$[\sigma_1^\xi]\in
R_d^\xi(p)/R_d(p)$ and $[\sigma_2^\xi]\in
R_d^\xi(p)/
\left\langle R_d(p),\sigma_1^\xi\right\rangle
$ are intrinsic, although their representatives are not. Changing coordinates or local generators replaces them by
$\sigma_1' \equiv u_1\sigma_1^\xi \mod R_d(p)$ and $\sigma_2' \equiv u_2\sigma_2^\xi+c\sigma_1^\xi \mod R_d(p)$, 
where $u_1,u_2$ are nowhere-vanishing regular functions. Hence their common span is unchanged and varies
regularly. Taking $a=0$ in the first-order computation gives
$$
\sigma_1^\xi\equiv wzf_x\mod R_d(p).
$$
No explicit representative of $\sigma_2^\xi$ is needed.
\end{remark}

\begin{proposition}\label{limd}
With $\Gamma$ as above, we have $\mu(B_d)\subset\Gamma$.
Moreover,
every element $u\in R_d^\xi(p)$ satisfies $\nu_p(u)\geq d-2$.
\end{proposition}

\begin{proof}
By \Cref{lem:two-step-limit-d}, the saturated family of restrictions defines,
locally along $E|_{C\smallsetminus\Gamma}$, a subbundle of rank $n_d$ of the
trivial bundle with fibre $S_d$. It therefore gives a regular extension of
$\widetilde\rho_d$, whose fibre at $\xi$ is
$$
R_d^\xi(p)
=
\left\langle
R_d(p),\sigma_1^\xi,\sigma_2^\xi
\right\rangle.
$$
Hence no point of $\mu^{-1}(C\smallsetminus\Gamma)$ belongs to $B_d$, and
$\mu(B_d)\subset\Gamma$.

Finally, every element of $R_d(p)$ has evaluation at least $d$, while
$\sigma_1^\xi$ and $\sigma_2^\xi$ have evaluations $d-1$ and $d-2$,
respectively.
\end{proof}

\section{The differential of the cone map}\label{sec: conemap}

\subsection*{Differential of the cone map}
 We now turn to the study of the differential of $\rho_d$.
We first give a standard definition.

\begin{definition} \label{defn:VmD}
Let $D$ be any effective divisor on $C$.
We set 
$$V(-D):= \{s\in V: (s)\geq D\},$$ 
i.e. the subspace of the section vanishing at $D$. Accordingly, for any $W\subseteq V$, we set 
$W(-D):=V(-D)\cap W.$ Note that $V(-p)=W(p)$ for any $p\in \pr^3$.

\end{definition}

We will now use an explicit setting analogous to \Cref{setting}, but for $p \in U$.

\begin{setting}\label{setting2}
$\,$
\begin{enumerate} 
\item $p=(0:0:0:1)$, a fixed point in $U$; 
\item $\ell=\{y=z=0\}$ the implicit equations of the line along which we deform;
\item $p_t=(-t:0:0:1)$, the corresponding linear deformation of $p$; 
\item $W(p_t)=\langle x+tw,y,z\rangle$ the corresponding deformation of $W(p)=\langle x,y,z \rangle$.
\end{enumerate}
\end{setting}

\begin{theorem}\label{limb}
The differential of $\rho_d$ is injective over $U$.
\end{theorem}
\begin{proof}
Let $p\in U$. We make the coordinate choice as in \Cref{setting2}, and consider the tangent vector $X$ as in \Cref{diff}. We have
$$
d_p\rho(X)\in\mbox{Hom}(R_d(p),S_d/R_d(p))
=T_{\rho(p)}G(n_d,S_d).
$$
Assume $d_p\rho(X)=0$ and, by contradiction, suppose that $X\neq0$. Let us consider $A(x,y,z)\in W(p)_d$, so that $A(x+tw,y,z)\in W(p_t)_d$. By taking the Taylor expansion as in Subsection \ref{subs:coordinates} we get \begin{equation}\label{eq: questa} A(x+tw,y,z)=A(x,y,z)+twX(A) \mod t^2. \end{equation}

Since $X\neq0$, for any $G(x,y,z)$ of degree $d-1$ we can find $A$ such that $X(A)=G$. In particular, from \Cref{eq: questa} we have $wG\in R_d(p)$ and therefore
$\phi_d(wW(p){d-1})\subseteq R_d(p)$. This means that
\begin{equation} \label{eq: questa2}
wW(p){d-1}\subseteq W(p)_d+I(d).
\end{equation}
We will show that this inclusion cannot happen.

Let $r_i$, for $i=1,2$, be two distinct bisecant lines through $p$, with $r_i\cap C=p_i+q_i$.
We can choose a basis
$W(p)=\langle x,y,z\rangle$ such that
$$z\in W(-(p_1+p_2+q_1+q_2)), \quad x\in W(-p_1-q_1), \quad y\in W(-p_2-q_2)$$
$$p_1, q_1 \notin \{y=0\}, \qquad p_2, q_2 \notin \{x=0\}.$$
We can also choose the coordinate $w$ such that
$$
w\in V(-(p_1+p_2)),
\qquad
q_1, q_2 \notin \{w=0\}.
$$
In other words, the points are in the following configuration:

\begin{figure}[h!]
\centering
\begin{tikzpicture}[
    scale=1.15,
    line cap=round,
    line join=round,
    point/.style={circle, fill=black, inner sep=1.6pt}
]

\coordinate (p)  at (0,0);
\coordinate (p1) at (2,1.25);
\coordinate (q1) at (3.2,2);
\coordinate (p2) at (1.35,-1.55);
\coordinate (q2) at (2.25,-2.58);

\draw[blue!70!black, very thick]
     (4.05,0.85)
    .. controls (3.3,1) and (3.55,1.50) .. (q1)
    .. controls (2.65,2.65) and (1.35,2.35) .. (p1)
    .. controls (2.65,0.45) and (2.55,-0.65) .. (p2)
    .. controls (0.45,-2.15) and (0.70,-3.15) .. (q2)
    .. controls (2.95,-2.40) and (3.65,-2.55) .. (3.65,-3.30);

\node[blue!70!black, right] at (2.55,-0.35) {$C$};

\draw[green!55!black, thick, dashed]
    ($(p1)!-0.55!(p2)$)
    --
    ($(p1)!1.55!(p2)$);

\node[green!45!black, left, align=right]
    at ($(p1)!0.45!(p2)$)
    {$\ell_w$};

\draw[magenta!80!black, thick]
    ($(p)!-0.75!(q1)$)
    --
    ($(p)!1.25!(q1)$);

\node[magenta!80!black, above left]
    at ($(p)!-0.40!(q1)$)
    {$r_1$};

\draw[orange!85!black, thick]
    ($(p)!-0.65!(q2)$)
    --
    ($(p)!1.25!(q2)$);

\node[orange!85!black, right]
    at ($(p)!-0.40!(q2)$)
    {$r_2$};

\foreach \P in {p,p1,q1,p2,q2}
    \node[point] at (\P) {};

\node[above, yshift=1mm]  at (p)  {$p$};
\node[right, xshift=2mm] at (p1) {$p_1$};
\node[above right, xshift=-2mm, yshift=1mm] at (q1) {$q_1$};
\node[left, xshift=-2mm]  at (p2) {$p_2$};
\node[above right] at (q2) {$q_2$};

\end{tikzpicture}
\caption{Representation of the configuration on the plane
$\{z=0\}$. Here $\ell_w=\{z=w=0\}$.}
\end{figure}

By \eqref{eq: questa2}, we have
$wx^{d-1}=s(x,y,z)+g(x,y,z,w)$,
where $s\in\sym^d\langle x,y,z\rangle$ and $g\in I(d)$.

Evaluating at $p_1$, we have $w(p_1)=x(p_1)=z(p_1)=0$ and $y(p_1)\neq0$. Since $g(p_1)=0$, it follows that $s$ does not contain the monomial $y^d$. Similarly, evaluating at $p_2$, we see that $s$ does not contain the monomial $x^d$. Hence
$$
s(x,y,z)=zs_1(x,y,z)+xys_2(x,y,z).
$$
In particular, $s(q_2)=0$, since $y(q_2)=z(q_2)=0$. Since also $g(q_2)=0$, evaluating the equality at $q_2$ gives
$$
w(q_2)x(q_2)^{d-1}=0,
$$
contradicting $w(q_2)\neq0$ and $x(q_2)\neq0$. Therefore $X=0$.
\end{proof}

The key point of the above  argument is that for $p\in U$ 
$$wW(p)_{d-1}\not \subseteq (W(p)_d+I(d)).$$
An equivalent way of phrasing this result is by stating that $\dim(\Gamma(p))>0$, where we set 
$$ \Gamma(p):= {\frac {w  W(p)_{d-1}}{(W(p)_d+I(d))\cap w W(p)_{d-1}}}.$$
By assuming $p$ to be general and $d$ to be large enough, we can prove a 
stronger result, which will be needed in \Cref{sec: conic divisors}.
\begin{proposition} \label{d(imp)} 
With the notations above, let $p\in U$ be a general point.
Then we have:
\begin{enumerate}
\item if $d>3$ then $dim(\Gamma(p))>1,$ 
\item if $d>4$ or $d =4$ and $g=0$ then $dim(\Gamma(p))>2.$  
\end{enumerate}
In case (2) we find a base of 
$W(p)=\langle x_1,x_2,x_3\rangle \subset \langle x_1,x_2,x_3,x_4 \rangle=V$ such that
$$x_4x_1^{d-1}, x_4x_2^{d-1}, x_4x_3^{d-1}$$ are independent$\mod (I(d)+W(p)_d).$
\end{proposition}
\begin{proof} Assume $d>3$; then for $p\in U$ general there are two distinct bisecant lines $r_i$, for $i=1,2$, that intersect $C$ in $p_i+q_i$. Then we can write $W:=W(p)=\langle x,y,z\rangle$ such that
$z\in W(-(p_1+p_2+q_1+q_2))$, $x\in W(-p_1-q_1)$ and $y\in W(-p_2-q_2)$. Finally, we take a general
$w\in V(-(p_1+p_2))$ such that $w(p),w(q_1),w(q_2)\neq0$.

We will show that $wx^{d-1}$ and $wy^{d-1}$ are independent in
$$
\frac {wW_{d-1}}{(I(d)+W_d)\cap w W_{d-1}},
$$
and this will imply $\dim(\Gamma(p))>1$. Assume by contradiction that there exist $a,b\in\bC$, not both zero, such that
$$
w(ax^{d-1}+by^{d-1})=g(x,y,z)\mod I(d),
$$
where $g\in W_d$. Since the left-hand side vanishes at $p_1$ and $p_2$, the polynomial $g$ cannot contain the terms $y^d$ and $x^d$, respectively. Hence we can write
$$
g(x,y,z)=zs(x,y,z)+xyt(x,y).
$$
Therefore $g$ vanishes at $q_1$ and $q_2$. Evaluating the relation at these points, and using $w(q_1),w(q_2)\neq0$, gives respectively $b=0$ and $a=0$, a contradiction.
\smallskip

Now we consider the case $d>4$ or $d=4$ and $g=0$. For a general point $p\in U$, the image curve
$C_p\subset\bP(W(p)^*)$
has $N=\frac{(d-1)(d-2)}2-g$
nodes. If $d=4$ and $g=0$, then $N=3$. Now by Castelnuovo's bound \cite[Chap. 2, Sec 3]{GH} $$N=\frac{(d-1)(d-2)}{2}-g\geq \frac{(d-1)(d-2)}{2}-\frac{(d-1)(d-3)}{4}=\frac{(d-1)^2}{4}.$$
These nodes cannot all lie on the same line in $\bP(W(p)^\ast)$. Indeed, the corresponding plane in $\bP^3$ would contain all the associated bisecant lines and would therefore intersect $C$ in at least $2N$ points. This is impossible, since $2N=6>d$ when $d=4$ and $g=0$, while
$$
2N\geq2\frac{(d-1)^2}{4}>d
$$
when $d>4$. We may therefore choose three non-collinear nodes, corresponding to three bisecant lines $r_i$, for $i=1,2,3$, through $p$ and not contained in the same plane.
\end{proof}

\section{Structural results on conic linear series}\label{sec: conic divisors}
We now develop a machinery to describe families of conic divisors on a non-degenerate smooth curve $C\subset \pr^3$ of degree $d$. 
In particular, we will prove  that the general degree $d^2$ divisor on $C$ is conic if and only if $d\leq 4$.

We begin by defining the space parametrizing cones in $\bP^3$ of degree $d$ which have vertex in $U$:
$$Y_d:=\bigcup_{p\in U} \bP (W(p)_d)\subset  \bP V_d.$$
Recall that $V=V_1=H^0(\pr^3,\cO_{\pr^3}(1))$ and for any $k\geq 1$, $V_k=H^0(\pr^3,\cO_{\pr^3}(k))$, and that we denoted by $W(p)$ the sections of $V$ vanishing at $p$, see \Cref{sec:preliminaries}.

Consider the tautological sequence of $ \cO_{U}(1)$ over $U$ 
\begin{equation}\label{eq: eulero}
0\la \cS \la V \otimes \cO_{U}\la \cO_{U}(1)\la 0.
\end{equation}

The stalk of $\cS$ over $p$ is $W(p)$. By taking the  the $d$-th symmetric product and the projectivization over $U$ of \eqref{eq: eulero}, we see that
 $Y_d$ is the image of the induced natural map $\tau_d$:

\[
\begin{tikzcd}
\bP_U(\sym^d(\cS)) \arrow[rr, "\tau_d"]\arrow[rd,twoheadrightarrow] & &  \bP(V_d)\\
& Y_d\arrow[ru,hookrightarrow] &
\end{tikzcd}
\]
Clearly $\tau_d$ is birational onto its image $Y_d$, and in particular we get:
\begin{equation} \label{eqn:dimU}
\dim Y_d = \dim \bP_U(\sym^d(\cS))=\binom{d+2}{2}+2.
\end{equation}
Moreover,  for any $p\in U$ we have that $W(p)_d\cap I(d)$ is a line in $W(p)_d$ and so it defines a line sub-bundle of $\sym^d(\cS)$ over $U$:
$$0\to \cL\to \sym^d(\cS).$$ 
Let $\cE:=\sym^d(\cS)/\cL$ be the quotient bundle. 
Observe that 
\begin{equation}\label{eq: dim}
\dim \bP_U(\cE)=\binom{d+2}{2}+1.
\end{equation}
\begin{remark}\label{rem: tautologico}
Note that $\cE$ could also be defined as the pullback via $\rho_d$ (\Cref{def: rho}) of the tautological bundle $\cP$ on the Grassmannian:
$\cE=\rho_d^*\cP$. Indeed, the fibre of $\rho_d^*\cP$ over $p\in U$ is 
$$\frac{W(p)_d}{W(p)_d\cap I(d)}\cong R_d(p) $$
by \Cref{rem: cone map}. We will use this point of view  in \Cref{def: phi''} to extend $\cE$ to (a suitable blow-up of) $\widetilde \pr^3$.
\end{remark}
We have an injection $j\colon U\to Y_d$ sending $p\mapsto [f_p]$, where $f_p$ is as usual an equation of the cone over $C$ with vertex in $p$. 
We set 
\begin{equation}
\cV:= Y_d\smallsetminus j(U) \subset  \bP(V_d).
\end{equation}
In other words, $\cV$ is the locus of cones in $\bP(V_d)$ which have vertices in  $U$ and are not cones over $C$.  
As $j(U)$ is closed in $Y_d$ we have that $\cV$ is an open subset of $Y_d$ and that $\dim \cV= \dim Y_d$. 

If $[f]\in \cV$, the restriction to $C$ gives  a well-defined morphism
\begin{equation}\label{Psi}
\Psi\colon \cV\la \bP(S_d),
\end{equation}
where we recall that $S_d:=\H^0(C,\cO_C(d))$.
The image of $\Psi$ consists of the conic divisors of degree $d^2$ over $C$ with vertex in $U$.

We remark that, if $[g]\in \bP(W(p)_d)\smallsetminus \{[f_p]\}$, then $\Psi$ is constant on the line 
generated by $[g]$ and $[f_p]$, since $f_p$ vanishes on $C$.
This implies that the rational map $\Psi\circ \tau_d\colon  \bP_U(\sym^d(\cS))\dasharrow \bP(S_d)$ factorizes to give a morphism $\Phi: \bP_U(\cE)\la \bP(S_d)$.
We summarise all this in the following commutative diagram:
\begin{equation}\label{Phi} 
\begin{tikzcd}
\bP_U(\sym^d(\cS)) \arrow[dd,twoheadrightarrow] \arrow[r,twoheadrightarrow] \arrow[rr,bend left=30,"\tau_d"] \arrow[rrd,bend right=50,dashrightarrow] & Y_d \arrow[r,hookrightarrow] & \bP(V_d)  \\
& \cV \arrow[r,"\Psi"] \arrow[u,hookrightarrow] & \bP(S_d)\\
\bP_U(\cE) \arrow[rru,bend right=30, "\Phi"]
\end{tikzcd}
\end{equation}

We want to understand when $\Phi$ and $\Psi$ are dominant. To this aim, we study the differential of $\Phi$.
\begin{theorem} \label{thm:diffmaxrank}
The differential of $\Phi$ is generically of maximal rank. 
More precisely, for the general point  $y\in \bP_U(\cE)$, if we let 
$$m:=\dim(\ker(d\Phi_y)),$$
we have that:
\begin{enumerate}
\item if $d\geq5$ or $d=4$ and $g=0,$ then $m=0$  ($d\Phi_y$ is generically injective);
\item  if $d=4$ and $g=1$ (i.e. $C$ is the elliptic normal curve), then $m=1;$
\item  if $d=3$ and $g=0$ (i.e. $C$ is the rational normal curve), then $m=2$.
\end{enumerate}
 \end{theorem}
\begin{proof}
Let  $\zeta\colon \bP_U(\cE)\to U$ be the fibre bundle map. 
Its differential gives the sequence
\[ 0\la T' \la T_{\bP_U(\cE)}\la T_U\la 0,\]
where $T'$ is the tangent along the fibres of $\zeta$ and $T_U= {T_{\bP^3}}_{|U}$.  
Fix a point $[g]\in \bP_U(\cE)$, call $p:=\eta([g]) \in U$ and assume $[g]\neq 0$, i.e. $g\not \in \langle f_p\rangle$.
We have 
$$T_{U,p} \cong \mbox{Hom} (W(p), V/W(p)).$$ 
Since $[g] \neq 0$, we have 
$$T'_{[g]}\cong W(p)_d/\langle g,f_p\rangle.$$
We remark that any fibre of $\bP_U(\cE)$ embeds in $\bP(S_d)$ via $\Phi$, so 
$$d\Phi_{|T'_{[g]}}\colon  T'_{[g]}\la T_{\bP(S_d),{[g]}}\cong S_d/\langle g\rangle$$ is injective. Recall that by \Cref{limnormali}, $R_d:=\phi_d(W(p)_d)$, 
we obtain 
$$d\Phi(T'_{[g]})=R_d/\langle g\rangle  \subset S_d/\langle g\rangle.$$

 We consider case (1):  $d>4$  or $d=4$ and $g=0$. 
 We have to find a point $[g]\in \bP_U(\cE)$ such that $d\Phi_{[g]}$ is injective: as we have already seen that $\Phi$ is injective on the fibres of $\zeta$, we have to find three ``horizontal'' tangent
 vectors $X_i$ in $[g]$ such that $d\Phi_{[g]}(X_i)\in S_d/\langle g\rangle$ are independent mod $d\Phi_{[g]}(T'_{[g]})=R_d/\langle g\rangle$.
 This is equivalent to the $d\Phi_{[g]}(X_i)$'s being independent in $S_d/R_d$. 
 
Arguing as in  Proposition \ref{d(imp)} (2), for the general point  $p\in U$ there are $3$ non-coplanar bisecants to $C$ passing through $p$, which we call $r_1,r_2,r_3$. 
We can choose coordinates 
 $x_1,x_2,x_3$ such that:
 \begin{itemize}
 \item $\langle r_j, r_k\rangle =\{x_i=0\}$, where $\{i,j,k\}=\{1,2,3\}.$ 
  \item $W(p)=\langle x_1,x_2,x_3\rangle$.
  \end{itemize} 
 Let us name the intersection points as follows: $r_i\cap C:=\{p_i, q_i\}$. 
Then we take
 $x_4$ such that $ V=\langle x_1,x_2,x_3,x_4\rangle$ and  $x_4(p_i)=0,\ 1\leq i\leq 3$.
 Note that with this choices the $p_i$'s are the three first coordinates points of $\pr^3$.

We fix a $[g]\in \bP_U(\cE)$ as follows:
$$g(x_1,x_2,x_3)=\frac{1}{d} (x_1^d+x_2^d+x_3^d).$$
Observe that $[g]$ is not zero, because $g(p_i)\not =0$ for $i=1,2,3$.

We define for any fixed triple of numbers $(a_1,a_2,a_3) \neq (0,0,0)$ the curve
$$\gamma(t):=[g(x_1+a_1tx_4,x_2+a_2tx_4,x_3+a_3tx_4)].
$$
Note that $\zeta(\gamma(t))$ is the point corresponding to the annihilator of the subspace $\langle x_1+atx_4,x_2+btx_4,x_3+ctx_4\rangle$.
Consider the tangent vector $X := \gamma'(0)$.
We have
$$d\Phi_{[g]} (X)=\left[a_1x_4x_1^{d-1}+a_2x_4x_2^{d-1}+a_3x_4x_3^{d-1}\right]\mod \langle g\rangle.$$
By \Cref{d(imp)} we have that $$x_4x_1^{d-1}, x_4x_2^{d-1}, x_4x_3^{d-1}$$ are independent 
modulo $(I(d)+W_d)$. 
As a consequence, their classes are independent in $S_d \mod R_d$. By a suitable choice of the $a_i$'s we can find three independent tangent vectors $X_i$ such that $d\Phi_{[g]}(X_i) = x_4x_i^{d-1}$. Hence $\ker(d\Phi_{[g]}))=0$.

The other cases can be  proved similarly by using \Cref{d(imp)} for the case $d=4$ and $g=1$, and \Cref{limb} for $d=3$. This concludes the proof. 
 \end{proof}

  \begin{corollary}\label{cor: dominant}
  The maps $\Phi$ and $\Psi$ are dominant $\iff$ $d\leq 4$.
  \end{corollary}
 \begin{proof}
 It is clear from Diagram \eqref{Phi},   $\Psi$ is dominant if and only if $\Phi$ is. 
 
Let us now consider the case when $\Phi$ has generically injective differential: using \Cref{thm:diffmaxrank} this is equivalent to $d\geq 5$ or $d=4$ and $g=0$.
In this case, $\Phi$ is dominant if and only if $\dim(\bP_U(\cE))\geq  \dim(\bP(S_d))$.
On one hand, 
\[\dim(\bP_U(\cE)) = \dim Y_d -1=\dim \cV-1=\binom{d+2}{2}+1,\] 
where the last equality is from \Cref{eqn:dimU}.
On the other hand, $\dim(\bP(S_d))=h^0(C,\cO_C(d))-1=d^2-g$. It follows that
$$\dim(\bP_U(\cE))\geq  \dim(\bP(S_d)) \iff g\geq d^2-\binom{d+2}{2}-1=\frac{d^2-3d-4}{2}.$$ 
From Castelnuovo's bound we have:
$g\leq \frac{(d-2)^2}{4}$. 
Putting these inequalities together we get that  necessarily $d\leq 4$; so we are in the case $d=4$, $g=0$. 
In this case $\dim(\bP_U(\cE))=16 =  \dim(\bP(S_d)),$
so $\Phi$ is dominant in this case, as wanted.
 
The differential of $\Phi$ is surjective when $d\leq 4$, $g\geq 1$ if  we verify that 
 \[\dim \cV=\dim \bP(S_d)+\dim\ker\Phi\] in the three cases  of \Cref{thm:diffmaxrank}, showing
  that $\Phi$ is generically submersive.
  This is easily proved. 
 For instance, if $d=4$ and $g=1$ we have $\dim\bP_U(\cE) =16$,  
 $\dim \bP(S_d)=15$, and $\dim \ker d\Phi=1$.  
 \end{proof}
 \begin{remark}\label{rem: dim conic divisors}
 If $d\geq 5$ the conic divisors on $C$ form a subvariety of dimension $\dim(\bP_U(\cE)) =\binom{d+2}{2}+1$ inside $\bP(S_d)$.
 \end{remark}

\subsection*{Partial compactification}
We partially extend $\Phi$ to (a suitable blow-up of) $\widetilde{\pr^3}$, by using the results of Section \ref{sec: limits}.

Let us consider the extension of the cone map to its maximal regular open set $\cV_d$ $$\widetilde{\rho_d}\colon \cV_d \longrightarrow G(n_d,S_d).$$ 
Recall that we proved in \Cref{limd} that the locus $B_d:=E\smallsetminus \cV_d$ on the exceptional divisor $E$ is contained in the union of fibres over $\Gamma$, where 
\[
\Gamma:=(C\smallsetminus C')\cup 
\{p\in C' \mid t_p\cap C\neq \{p\}\}.
\]
By taking the pull-back of the tautological bundle over $G(n_d,S_d)$ as in \Cref{rem: tautologico}, we obtain $\cE'$.
Take the projectivization $\bP_{\cV_d}(\cE')$
and  extend  $\Phi$ to a map
$$\Phi' \colon \bP_{\cV_d}(\cE') \longrightarrow \bP (S_d).$$

Recalling that $\nu\colon \widetilde{\pr^3}\to \pr^3$ is the blow-up of $\pr^3$ along $C$, note that the complement of $\cV_d$ outside $E$ is contained in the finite set $\nu^{-1}(S\smallsetminus C)$. 
Let $\tau \colon \widehat \pr^3 \to G(n_d,S_d)$ be a suitable series of blow-ups that extends $\widetilde{\rho_d}$ to $\cV_d':=\widehat\pr^3 \smallsetminus B_d $:
\[
\widehat{\rho_d}\colon\cV_d'\longrightarrow  G(n_d,S_d).
\]
\begin{definition}\label{def: phi''}
Define  $\cE''$ to be the extension of $\cE'$ over $\cV_d'$ obtained via $\widehat{\rho_d}$,
 and $\Phi''$ to be the corresponding extension of $\Phi'$: 
\begin{equation}\label{eq: phi''}
    \Phi'' \colon \bP_{\cV_d'}(\cE'') \longrightarrow \bP (S_d).
\end{equation}
\end{definition}

\section{An application: some special pencils of plane quartics} \label{sec:application}

Take an elliptic curve $C$ with origin $O,$ define the  embedding  $j: C\to |4O|=\bP^3,$ so that $j(C)$ is the  elliptic normal curve and $O$ is an  inflection point of $j(C)$. 
We  identify  $C$ and $j(C)$.  
Recall that it is a complete intersection of two quadrics  $$C =  \{Q_1=Q_2=0\}.$$
\begin{remark} \label{rmk:nopoints}
Let us find the non-birational projection points $S$ and see that they are not contained in $C$.
Indeed $S$ consists of the vertices of the four singular quadrics in the pencil 
$$\lambda Q_1+\mu Q_2, \qquad (\lambda:\mu)\in \pr^1.$$ 
After a suitable coordinate choice, we have 
$$Q_1=x_1^2+x_2^2+x_3^2+x_4^2;\ Q_2=a_1x_1^2+a_2x_2^2+a_3x_3^2+a_4x_4^2,\ \
a_i\neq a_j, \iff i\neq j .$$
The singular quadrics are given by $(\lambda, \mu)\equiv (-a_i,1)$ and the vertices are the coordinate points. 
\end{remark}

We can prove in this setting the following result.
 \begin{theorem}\label{do}
 If $C\subset \pr^3$ be as above.
Then the map $\Phi''$ is proper and surjective.
\end{theorem}
\begin{proof}
We already proved in  \Cref{cor: dominant} that  for $d\leq 4$ the map $\Phi$ is dominant. 
We only need to observe that for an elliptic space curve of degree $4$ the set $\Gamma$ is empty. Indeed, by \Cref{rmk:nopoints} above,
$S\cap C=\emptyset$. Moreover, for every $p\in C$ one has $t_p\cap C=2p$ scheme-theoretically. Indeed, since $C=Q_1\cap Q_2$, both quadrics restrict to sections of degree $2$ on $t_p$ vanishing at least twice at $p$. If the intersection were larger than $2p$, both restrictions would vanish identically, so $t_p\subset C$, a contradiction. In particular, $\pi_p(p)$ is smooth on $C_p$. Hence $\Gamma=\emptyset$.
\end{proof}

\begin{remark}\label{rem: order}
Let  $p$ be any point on  $E\smallsetminus \nu^{-1}(S)\subset \widetilde \pr^3$. 
By \Cref{limd}, any divisor in the image via $\Phi'$ of the fibre of $\cE'$ has vanishing on $p$ of order at least $d-2$. 
\end{remark}

  \begin{corollary}\label{ord}
 Let  $p\in C$. If  $D\in \Phi'\big( \bP(\cE'_{|\bP(N_p)})\big)$, i.e. $D$ belongs to $\pr(R_d^\xi(p))$, for some $\xi\in E$ such that $\eta(\xi)=p$. 
 Then $D=(d-2)p+R$ where $R$ is effective. 
  \end{corollary}
  \begin{proof}
By applying part (2) of \Cref{limd} it follows that $D-(d-2)p$ is still effective. This concludes the proof.
  \end{proof}
  \begin{remark}
 For any point $p\in C'$ we let $\bP(N_p)\cong \bP^1$ be the projectivization of the normal bundle $N_p$. 
 We have that the restriction
  \begin{equation} \label{eqnEprimo}
{\cE'}_{|\bP(N_p)}\to \bP(N_p)
\end{equation}
is a vector bundle.
Making the point $p$ vary in $C'$ in \Cref{eqnEprimo} we get a natural map coming from the bundle projection
\begin{equation} \label{eqn:mu}
   \eta\colon \bP_{U'}(\cE')\smallsetminus \bP_U(\cE)\longrightarrow C.
\end{equation}
 \end{remark}

We have $\dim H^0(C,\cO_C(4))=16 $ and so $\bP(S_4)\cong \bP^{15}$.
We recall that  $\dim\cV=17$ (see \Cref{eqn:dimU}), and that the maps $\Psi\colon\cV\to \bP(S_4)$  and $\Phi\colon\bP_U(\cE)\to \bP(S_4)$ are both dominant.

As proved in \Cref{do}, 
$\Phi''\colon \bP_{\widehat{\pr^3}}(\cE'')\longrightarrow \bP(S_4)$ is surjective.
Thus we  know that a general divisor in $\pr(S_4)$ is a conic divisor and that for any effective divisor  $D \in \pr(S_4)$ the inverse image ${\Phi''}^{-1}(D)$ is non empty. 

Now we prove that \emph{any} divisor $D$ in $\pr(S_4)$ is a conic divisor, that is ${\Phi''}^{-1}(D)\cap \pr_U(\cE)$ is not empty.
Before proving this result, let us establish this simple lemma:

\begin{lemma} \label{lemma:dimfibre}
Let $X,Y,Z$ be smooth projective varieties, and $\varphi\colon X\to Y$ and $\eta\colon X\to Z$ be proper surjective morphisms such that for a general $y\in Y$ we have
\begin{equation} \label{eqn:dimensionlem}
    \dim(\eta(\varphi^{-1}(y))=\dim \varphi^{-1}(y)>0.
\end{equation}
Then $\dim(\eta(\varphi^{-1}(y)))\geq 1$ holds for any $y \in Y$.
\end{lemma}
\begin{proof}
Let $W:=(\varphi,\eta)(X)\subseteq Y\times Z$. Since $\varphi$ is proper and $Z$ is separated, the morphism $(\varphi,\eta)\colon X\to Y\times Z$ is proper. Hence $W$ is closed and the induced projection $W\to Y$ is proper. Its fibre over $y\in Y$ is precisely
$W_y:=\eta(\varphi^{-1}(y))$.
By upper semicontinuity of fibre dimension, \cite[Chapter II, Exercise 3.22(d)]{Hartshorne}, the locus
$$
\{y\in Y\mid \dim W_y\geq 1\}
$$
is closed. By assumption it contains a dense open subset of $Y$, hence it is all of $Y$. So $\dim W_y\geq 1$
for every $y\in Y$.
\end{proof}

\begin{theorem}\label{teo: contro}  
Fix an effective divisor $D \in \pr(S_4)$ on $C$. 
There is a $1$-dimensional variety $B\subset \bP_U(\cE)$ such that:
\begin{enumerate}
\item[(i)] the fibre projection of $B$ over $U$ is again $1$-dimensional;  
\item[(ii)] by calling $\widetilde{K}_t$ the cone in $\pr^3$ corresponding to $t\in B$, we have that  $\Phi''(\widetilde{K}_t)=\Phi(\widetilde{K}_t)= D$ for all $t\in B$.
\end{enumerate}
\end{theorem}
\begin{proof}
Set $G:=(\Phi'')^{-1}(D)\subset \bP_{\widehat{\pr^3}}(\cE'')$. The set $G$ is non-empty because $\Phi''$  is surjective. 
Since $\dim \bP_{\widehat{\pr^3}}(\cE'') = \dim \bP(S_4)+1$ by \Cref{eq: dim}, we have that $\dim G\geq 1$.
We apply \Cref{lemma:dimfibre} with $X:=\bP_{\widehat{\pr^3}}(\cE'')$, $Y:=\pr(S_4)$ and $Z:=\pr^3$. The map $\varphi$ is $\Phi''$ and the map $\eta$ is the restriction of the fibre projection of $X$ composed with the blow down to $\pr^3$. 
It is clear that for $p \in U$ the hypothesis in \Cref{eqn:dimensionlem} holds, so the dimension of $\eta(G)$ is greater or equal than $1$. 
We want to prove that $\eta(G) \cap U \neq \emptyset$. 
Assume by contradiction that $\eta(G) \subset S \cup C$. Since $\dim \eta(G) \geq 1$, necessarily $C$ must be contained in $\eta(G)$. 
However, by \Cref{ord} $\eta(G)$ must be contained in the support of $D$, a contradiction. So, there is at least a $1$-dimensional complete subvariety $\overline B$ of $G$ such that $\eta(\overline B) \cap U \neq \emptyset$. The corresponding $B=\overline B\cap  \bP_U(\cE)$  is the family we want. Properties $(i)$ and $(ii)$ are clear from the construction.
\end{proof}
In the above theorem we proved the existence of  cones in $\pr_U(\cE)$, i.e. cones with vertex in the open $U$ cutting on $C$  any divisor $D \in \pr(S_4)$. 
For what follows we need to study elements in the complement $\bP_{U'}(\cE')\smallsetminus \bP_U(\cE)$, which correspond to a specific divisor $D$, under the map $\Phi''$.

\begin{lemma}\label{lemma: conovero}  
Fix $q\in C$ such that $16q\sim 16O$ and set $D:=16q$. Let $\pi_q\colon C \to \pr^2$ be the linear projection and $G:=(\Phi'')^{-1}(D)\subset \bP_{\widehat{\pr^3}}(\cE'')$.
Consider the set 
$$\Gamma := G \cap ( \bP_{U'}(\cE')\smallsetminus \bP_U(\cE))$$
whose fibres parametrize  the elements in $\pr(R^\ell_p)$, with $p\in C$, $\ell \neq t_p$ containing $p$, which are $16q$ as divisors on $C$. 
Then $\Gamma$ is a rational curve and every element of $\Gamma$ is the class of a  cone over a fixed quartic $\cQ \subset\pr(W(q)^\ast)\cong  \pr^2$ intersecting $C$ in $16q$.
\end{lemma} 
\begin{proof} 
As in the proof of \Cref{teo: contro}, $G$ is non-empty and $\dim G\geq 1$.
Recall that  $\eta \colon ( \bP_{U'}(\cE')\smallsetminus \bP_U(\cE)) \to C$ is  the restriction of the  fibre projection.
Since $D$ has support only on $q$, from \Cref{ord} we get $\Gamma \subseteq \eta^{-1}(q)$. Indeed, if $K$ is an element of $\Gamma$, the divisor $\Phi'(K)$  belongs to $\pr(R^\ell_4(p))$ for some $p$ in $C$ and some line $\ell$ with direction in $\pr(N_p)$, so by \Cref{ord}  is of the form $2p + R$, with $R$ effective. On the other hand by assumption $\Phi'(K)=16q$, so necessarily $p=q$.

The linear system $|4O-q|$
gives the linear projection $\pi_q\colon C\to \pr(W(q)^\ast)\cong \pr^2$, with image a smooth cubic $C_q=\pi_q(C)$. 
The map $\pi_q\colon C\to C_q$ is an isomorphism.
Let $\overline O$ be the image of $O$ and $\overline q$ be the image of $q$.
We have by assumption that $16\overline q\sim 16 \overline O$ on $C_q$, and hence 
$$4(4\overline O-\overline q) \sim 12\overline q.$$ 
This linear equivalence implies that there exists a quartic curve $\cQ \subset \pr(W(q)^\ast)$ intersecting $C_q$ in $12 \overline q$. 

Now observe that the cone $K_\cQ\subset\pr^3$ corresponding to $\cQ$ cuts the curve $C$ in $16q$. 
Indeed, the elements in $\pr(W(q)_4)$ have a vanishing of order $4$ in $q$. 
So the class of $K_\cQ$ is an element of $G\smallsetminus \pr_U(\cE)$, and being a cone, its class also is an element of the linear system $\pr(R_4(q))$. 
Now recall from \cref{limd} that $R_4(q)$ is  contained in all the $R_4^\ell(q)$'s for all lines $\ell$ with directions parametrised by the exceptional divisor $\bP(N_q)$.
Hence the class of $K_\cQ$ defines a (constant) section of  
$$\bP_{\widetilde{\pr^3}}(\cE'')_{|\bP(N_q)}\to \bP(N_q)\cong \bP^1,$$ 
and therefore a rational curve $\Gamma\subseteq G\cap Z.$ 
Now recall that the restriction of $\Phi''$ to the fibre $\pr (\cE''\otimes \mathbb C(q) )\cong \pr (R^\ell_{4}(q'))$ to $\pr(S_4)$ is injective, and so we have that $\Gamma=G\cap Z$. This concludes the proof.
\end{proof}

We are now able to prove the main result of this section, which provides a degree $4$ analogue of \cite[Theorem 1.3]{MPS}.

\begin{theorem}\label{thm: pencil}
There exists a pencil of quartics in $\bP^2$ such that
 \begin{enumerate} 
  \item the base locus is set-theoretically one point;
 \item all the quartics of the pencil are irreducible;
\item the general element of the pencil is smooth;
 \item the pencil is non-isotrivial.
   \end{enumerate}
 \end{theorem}
\begin{proof}

\noindent{\it Proof of (1)}.
 As above, we take an elliptic curve $C$ with origin $O$, and  define the  embedding 
 \[
j\colon C\longrightarrow |4O|=\bP^3,
 \]
 so that $j(C)$ is the  elliptic normal curve and $O$ is an inflection point of $j(C)$. 
 For simplicity write $C=j(C)$.  
 Fix a point $q\in C$
such that $16q\sim 16O$ but $8q\not \sim 8O$;
so, $q$ is a torsion point of order exactly $16$.
By applying \Cref{teo: contro} to $D=16q$, we find a $1$-dimensional family of quartic cones, which we call  $\widetilde K_t$, with $t\in B\subset \pr_U(\cE) $, such that $\widetilde K_t\cdot C=16q, $ for any $t\in B$. Let
$\widetilde K_t=\{k_t=0\}$, so that $\Phi([k_t])=16q$. 
Fix $\widetilde K=\{k=0\}$ a general one of these cones, and let $p\in U$ be its vertex. 
Let $\Pi_p\colon \bP^3\smallsetminus\{p\}\to \bP^2$ be the projection. We have that $\Pi_p(\widetilde K\smallsetminus\{p\})=K_p$ is a quartic and so is the projection $\Pi_p (C)=C_p$. 
 
The pencil we want is defined by $\lambda k+\mu f_p=0$ in $\bP^2$. It satisfies (1) by construction, because set-theoretically $C_p\cap K_p=\{\Pi_p(q)\}$. Let $T_{\lambda,\mu}\sim\lambda K_p+\mu C_p$ be an element of this pencil. 

\noindent{\it Proof of (2)}. If this were not the case, 
then some $T_{\lambda,\mu}$ would be either the union of a line and a cubic or of two conics. In the first case we would have
$4O\sim4q$ and in the second case $8O\sim 8q$, both contradicting our assumption on $q$.

\noindent{\it Proof of (3)}.  This is the most delicate rgument of our proofAssume by contradiction that all the pencils constructed above have only singular members. 
Let us fix one of these pencils. It is generated by the projection from $p\in U$ of $C$ and by the projection of a cone $K_{t'}$, with $t'\in B$, with vertex $p$. It is convenient to fix in this case a plane $\Pi\cong \pr^2$ in $\pr^3$ not containing $p$ nor $q$, and consider the projections over this fixed frame.

By Bertini's theorem the only common singularity of the elements of the pencil must be at the point $\pi_p(q)=\overline q$. 
Then $C_p$ is singular in $\overline q$ only if the line $\ell$ joining $p$ and $q$ is the tangent line to $C$ at $q$.
By construction, for all $t\in B$, all the cones $\widetilde K_{ t}$ cut $C$ in $16q$ and by the contradiction assumption they are all singular in $\overline q$. This means that  $\ell$ is also contained in all the other cones $\widetilde K_{ t}$, for any $t\in B$.
In other words, the fibre projection of $B$ over $U$ is contained in $U\cap \ell$.
Since $B$ corresponds to a $1$-dimensional family of quartic cones cutting $C$ on $16q$ and all the pencils associated to any of these cones have all  members singular, we conclude that the vertices of these cones must all lie on $\ell$.

We  now specialise the projection point $p$ to $q$ along $\ell$. 
Call $t=0$ the point of the closure of $B$ corresponding to $q$.
Consider the plane  pencil  $\cP$  in $\Pi$ generated by the curves obtained by projecting  $\widetilde K_0^\ell$ and the limit cone $\widetilde C_0^\ell$ from $q$.

Denote by $\pi_q\colon C\to \Pi$ the extension of the projection from $q$, and $C_q\subset \Pi$ the image curve, which is a smooth plane cubic isomorphic to $C$. 
Observe that $C_q$ passes through $\overline{q}$ and has a simple tangent in this point. 
Indeed, the osculating plane in $\pr^3$ to $C$ at $q$ intersects $C$ in the divisor $3q+q'$, with $q'\not = q$, because of the hypothesis $4q \not \sim 4O$. 

By \Cref{lem:osculating}, the projection of $\widetilde C_0^\ell$ is the reducible quartic $C_q\cup L$,  where $L$ is the tangent line to $C_q$ at $\overline q$ (the projection of the osculating plane). 
The intersection divisor with $C_q$ is thus $2\overline{q}+\pi_q(q')$.

Since we have taken a limit of $K_t$ along $B$ going to a point over $C$, we have that $K_0^\ell$ necessarily lies in the rational curve $\Gamma$ that we found in \Cref{lemma: conovero}.
The pencil $\cP$ is thus generated by the quartic $\cQ$ of \Cref{lemma: conovero} and by $C_q \cup L$. 
The general element of $\cP$ is necessarily reduced because $C_q \cup L$ is. 
Moreover, \Cref{lemma: conovero} shows that the base locus of $\cP$ is set-theoretically $\overline q$.

We know that $D \cdot (C_q+L)=16\overline q$, and as a consequence $D\cdot C_q=12\overline q$ and $D \cdot L=4\overline q$. 
Therefore $D \cdot 3L=12\overline q$ and so there is a pencil of cubics osculating $D$ at $\overline q$ generated by $C_q$ and $3L$. 
This pencil covers $\bP^2$, so there exists a member of the pencil 
passing through $\overline q$ and another point $q'$ of $D$. This is a contradiction.

\noindent{\it Proof of (4)}.  To show that the pencil constructed is non-isotrivial, we remark that it contains a reduced member whose geometric genus is $1$ (the image of $C$ itself). 
The stable reduction of this member  necessarily is the same singular curve, while by point (3) above the general member is smooth.
\end{proof}

\begin{remark}\label{ni}
Observe that for any degree $d$ smooth curve in $\pr^3$ there does {\em not} exist an irreducible curve $C'\subset \pr^3$ of degree $d'$ such that the intersection of $C$ with $C'$ consists of $dd'$ points counted with multiplicity, see \cite{Diaz, MR2174485, MR3320214}. 
Instead, under the assumptions of Theorem \ref{thm: pencil}, we have proved  that there exists a cone over a smooth plane quartic that intersects $C$ in $p$ with multiplicity $16$.
\end{remark}

\begin{remark}
The properties of the members of the pencils studied in this section are related to the higher Weierstrass points, as studied for instance in \cite{KE}. 
There, it is discussed an explicit construction of a smooth plane quartic curve $\cQ$ that intersects a plane cubic at a point $q$ such that $4(4O-q) \sim 12q$ but $8O$ is not linearly equivalent to $8q$.
\end{remark}

\section{Conclusions and open problems}
The technique of using cones to study specific divisors on curves seems very promising. We conclude our work with a list of open problems highlighting potential directions for future research.

\begin{problem}
Study the map $\Phi$ for degree greater than $5$ and for curves embedded in $\pr^n$, with $n\geq 4$. 
As observed in \Cref{rem: dim conic divisors}, the closure of the image of $\Phi$ is a subvariety of $\pr(S_d)$ of codimension $d^2-g-\left(\binom{d+2}{2}+1\right)=\frac{(d-1)(d-4)}{2}-g$.
It would be interesting to find geometric conditions under which a divisor belongs to  this subvariety.
\end{problem}

A precise description of the limits of the conic linear series in higher degree could be interesting:

\begin{problem}
    Find results analogous to \Cref{cor:evaluation-d-1} and \Cref{limd} for cones of degree $d+1$ and higher.
\end{problem}

As for the results in \Cref{sec:application}, it would be interesting to generalize the setting to other curves which are not elliptic, for instance:

\begin{problem}
Construct a similar framework as in \Cref{sec:application}, starting with $C$ a (possibly non smooth) rational curve in $\pr^3$. The first problem is to check if we have in this case a result similar to  \Cref{do}.
\end{problem}

\begin{problem}\label{cir1}[Suggested by Ciro Ciliberto]
Consider the linear series on a curve cut out by a surface of degree $d$ having a multiple point whose multiplicity is not necessarily $d$.
\end{problem}

\begin{problem}\label{cir2}[Suggested by Ciro Ciliberto]
 One could then
consider ``conic series'' for varieties of dimension greater than
$1$, such as surfaces.
\end{problem}

\begin{problem}
The pencils constructed in  \Cref{sec:application} are degree $4$ integrable foliations on the projective plane minus one point. We wonder if some of our results can be related to the study of foliations as done for instance in \cite{AZ} by C. Alc\'antara and A. Zamora.
\end{problem}

Finally, in the recent paper \cite{weichen}, W. Chen uses the pencil of cubics with one base point studied in \cite{MPS} to give an example related to the hyperbolicity of the Hirzebruch surface $\mathbb{F}_1$ minus a curve $B$ such that $(\mathbb{F}_1,B)$ is of log general type. This is related to the Lang conjecture for pairs of log general type $(\pr^2, B)$ where $B$ is a curve with at least three components, as treated for example in \cite{C-T}.

\begin{problem}
Investigate whether the generalization of \cite{MPS} to quartics given in \Cref{thm: pencil} can be applied to problems related to the hyperbolicity of the complement of a curve in a rational surface.
\end{problem}

\bigskip

\noindent Riccardo Moschetti,\\Dipartimento di Matematica, Universit\`a di Pavia, Italy.\\
E-mail: \textsl {riccardo.moschetti@unipv.it}.
\medskip

\noindent Gian Pietro Pirola,\\Dipartimento di Matematica, Universit\`a di Pavia, Italy.\\
E-mail: \textsl {gianpietro.pirola@unipv.it}.
\medskip

\noindent Lidia Stoppino,\\Dipartimento di Matematica, Universit\`a di Pavia, Italy.\\
E-mail: \textsl {lidia.stoppino@unipv.it}.


\begin{thebibliography}{99}

\bibitem{AZ} C.R.~Alc\'antara, A.G.~Zamora, {\em  Some families of pencils with a unique base point and their associated foliations}, Qual. Theory Dyn. Syst. {\bf 23} (2024), no.~5, Paper No. 212, 18 pp.

\bibitem{ACGH}  E.~Arbarello, M.~Cornalba, P.A.~Griffiths, J.~Harris, {Geometry of algebraic curves. {V}ol. {I}}, {GMW}, {\bf 267}, {Springer-Verlag, New York}, {1985}.

\bibitem{BL} J.~Blanc, S.~Lamy, {\em Weak Fano threefolds}, Proc. LMS (3)  {\bf 105}  (2012), 1047--1075.

\bibitem{BoP} M.~Bolognesi, G.~P.~Pirola, {\em Osculating spaces and Diophantine equations (with an Appendix by Pietro Corvaja and Umberto Zannier)}, Math. Nachr. {\bf 284} (2011), no.~8-9, 960--972.

\bibitem{C-T} L.~Caporaso, A.~Turchet, {\em Hypertangency of plane curves and the algebraic exceptional set}, Proc. Lond. Math. Soc. (3) {\bf 130} (2025), no.~6, Paper No. e70063.

\bibitem{Cast} G.~Castelnuovo, {\em Sui multipli di una serie lineare di gruppi di punti appartenente ad una curva algebrica}, Rend. Circ. Mat. Palermo {\bf 7}, 89-110 (1893).

\bibitem{Cayley} A.~Cayley, {\em On reciprocants and differential invariants}, Quarterly Journal of Pure and Applied Mathematics, xxvI. (1893), 169--194, 289--307.

\bibitem{weichen} W.~Chen, {\em Algebraic Exceptional Set of a Three-Component Curve on Hirzebruch Surfaces}, Int. Math. Res. Not. I(2026), no. 11, Paper No. rnag113.

\bibitem{CCM} M.G.~Cifani, A.~Cuzzucoli, R.~Moschetti, {\em Monodromy of projections of hypersurfaces}, Ann. Mat. Pura Appl. {\bf 201}  (2022), 637--654.

\bibitem{CM} M.G.~Cifani, R.~Moschetti, {\em A note on non-uniform points for projections of hypersurfaces}, Annali dell'Università di Ferrara. {\bf 68}  (2022).

\bibitem{CF} C.~Ciliberto, F.~Flamini, {\em On the branch curve of a general projection of a surface to a plane}, Trans. Amer. Math. Soc. {\bf 363} (2011), no.~7, 3457--3471.

\bibitem{CMR} C.~Ciliberto, M.~Mella, F.~Russo, {\em Varieties with one apparent double point}, J. Algebraic Geom. {\bf 13} (2004), no.~3, 475--512.
 
\bibitem{Diaz}  S.~Diaz, {\em Space curves that intersect often}, Pacific J. Math., {\bf 123} (1986)  no.~2, 263--267.

\bibitem{Gat} R.~Gattazzo, {\em Points of type 9 on an elliptic cubic (Italian)}, Rend. Sem. Mat. Univ. Padova {\bf 61} (1979), 285--301. 

\bibitem{GH}  P.H.~Griffiths, J.~Harris, Principles of algebraic geometry, Pure Appl. Math. Wiley-Interscience, New York, 1978, xii+813 pp.

\bibitem{GLP} L.~Gruson, R.~Lazarsfeld, C.~Peskine, {\em On a theorem of Castelnuovo, and the equations defining space curves}, Invent. Math. {\bf 72} (1983), no.~3, 491--506.

\bibitem{Hartshorne} R. Hartshorne, {\em Algebraic geometry}, Graduate Texts in Mathematics, No. 52, Springer, New York-Heidelberg, 1977; MR0463157

\bibitem{MR3320214} R. Hartshorne, R.M.~Mir\'o-Roig, {\em On the intersection of ACM curves in $\mathbb{P}^3$}, J. Pure Appl. Algebra {\bf 219} (2015), no.~8, 3195--3213.

\bibitem{K} H.~Kaji, {\em On the tangentially degenerate curves}, J. London Math. Soc. {\bf 33} (1986), no. 3, 430--440.

\bibitem{KE} A.~Kamel, W.K.~Elshareef, {\em Weierstrass points of order three on smooth quartic curves}, J. Algebra Appl. {\bf 18} (2019), no. 1, 1950020, 21 pp.

 \bibitem{Kollar} J.~Koll\'ar, {\em Log K3 surfaces with irreducible boundary}, arXiv:2407.08051 [math.AG].

\bibitem{Laz1} R.~Lazarsfeld, {Positivity in algebraic geometry. {I}}, Ergebnisse der Mathematik und ihrer Grenzgebiete. 3. Folge. {48}, Springer-Verlag, Berlin, 2004, pp.xviii+387.

\bibitem{MPS} R.~Moschetti, G.P.~Pirola, L.~Stoppino, {\em Pencils of plane cubics with one base point}, Rend. Circ. Mat. Palermo (2) {\bf 74 } (2025), no.~1, Paper No. 60, 12 pp.

\bibitem{MR2174485} R.~M. Mir\'o-Roig and K. Ranestad, Intersection of ACM-curves in $\mathbb{P}^3$, Adv. Geom. {\bf 5} (2005), no.~4, 637--655.

\bibitem{PS}  G.P.~Pirola, E.~Schlesinger, {\em Monodromy of projective curves}, J. Algebraic Geom., {\bf 14} (2005), 623--642.

\bibitem{szip} L.~Szpiro,   Lectures on equations defining space curves, Tata Institute of Fundamental Research, Bombay; Springer-Verlag, Berlin-New York, 1979, v+81 pp.

\end{thebibliography}
\end{document}